\documentclass[12pt]{article}
\usepackage{dsfont}
\usepackage{mathrsfs,amssymb,amsmath}

\title{Commutative Stochastic Games}
\author{Xavier Venel}
\usepackage[pdftex]{graphicx}
\DeclareGraphicsExtensions{.jpg,.ps,.eps}
\usepackage{epstopdf}%
\newtheorem{definition}{Definition}[section]
\newtheorem{theorem}[definition]{Theorem}
\newtheorem{proposition}[definition]{Proposition}

\newtheorem{lemma}[definition]{Lemma}

\newtheorem{corollary}[definition]{Corollary}
\newtheorem{remark}[definition]{Remark}
\newtheorem{example}[definition]{Example}

\newcommand{\1}{\mathds{1}}

\newcommand{\Z}{\ensuremath{\mathds{Z}}}
\newcommand{\PP}{\ensuremath{\mathds{P}}}

\newcommand{\E}{\ensuremath{\mathds{E}}}

\newcommand{\N}{\ensuremath{\mathds{N}}}

\newcommand{\R}{\ensuremath{\mathds{R}}}

\renewcommand{\P}{\ensuremath{\mathds{P}}}
\renewcommand{\phi }{\varphi }
\newcommand{\T}{\ensuremath{\mathcal T}}

\author{Xavier Venel \footnote{School of Mathematical Science, Tel Aviv University, {\rm
email: xavier.venel@gmail.com}}
   }

\begin{document}
\maketitle
%

%We are interested in stochastic games with finite sets of actions where the transitions commute. It means that given the current state playing the action profile $a1$ followed by the action profile $a2$ whatever is the state at stage $1$ leads to the same distribution of state as playing first the action profile $a2$ and then $a1$. The Big Match and more generally absorbing games can be formulated in this model. When there is only one player and the transition mapping is deterministic, we show that the existence of a uniform value in pure strategies implies the existence of $0$-optimal strategies. In the framework of two-player stochastic games, we study one class of stochastic game with signals. When players observe past actions but not the state, we prove that the game has a uniform value. These games reduce to a specific class of zero-sum stochastic games on $\R^n$, where given a profile of action the transition is deterministic. We study this auxiliary game by induction on the number of "non cyclic" actions and by using the theorem of Mertens and Neyman (1981). The same proof extends to the non zero-sum case if we use the result of Vieille (2000).
%
\begin{abstract}
We are interested in the convergence of the value of $n$-stage games as $n$ goes to infinity and the existence of the uniform value in  stochastic games with a general set of states and finite sets of actions
where the transition is commutative. This means that playing an action
profile $a_1$ followed by an action profile $a_2$, leads to the same
distribution on states as playing first the action profile $a_2$ and
then $a_1$. For example, absorbing {games can} be reformulated as
commutative stochastic games.

%Many problems satisfy this assumption or can be reformulated in order to satisfy it such as absorbing games.

When there is only one player and the transition function is
deterministic, we show that the existence of a uniform value in pure
strategies implies the existence of $0$-optimal strategies. In the
framework of two-player stochastic games, we study a class of games where the set of states is $\R^m$ and
the transition is deterministic and $1$-Lipschitz for the
$L_1$-norm, and prove that these games have a uniform value. A similar proof shows
the existence of an equilibrium in the non zero-sum case.

These results remain true if one considers a general
model of finite repeated games, where the transition is commutative and the
players observe the past actions but not the state.

%Given one of these
%games, one can introduce an auxiliary stochastic game, deterministic
%and commutative, which satisfies the assumptions of the previous
%paragraph.
\end{abstract}
\normalsize

%\title{Commutative Stochastic Games}

%\bibliographystyle{MOR}
%\maketitle

%-------------------------------------------------------------------------------------------------------------------
%\citet{Grandell}
%\citep[see][]{Grandell}
% \newpage
% \tableofcontents

%\newpage
\section{Introduction}
A two-player zero-sum repeated game is a game played in discrete time. At each stage, the players independently take some decisions, which lead to an instantaneous payoff, a lottery on a new state, and a pair of signals. Each player receives one of the signals and the game proceeds to the next stage.

This model generalizes several models that have been studied in the literature. A Markov Decision Process (MDP) is a repeated game with a single player, called  a decision maker, who observes the state and remembers his actions. A  Partial Observation Markov Decision Processes (POMDP) is a repeated game with a single player who observes only a signal that depends on the state and his action. A stochastic game, introduced by Shapley \cite{Shapley_53}, is a repeated game where the players learn the state variable and past actions.\\

%If the players do no observe the state but the signal is public, i.e. players observe the same signal, the game is called a symmetric repeated game. In this paper, we will focus on MDPs, Stochastic games and on POMDPs and repeated games where the players observes past actions played but no information on the state. \\

%trivial signal, and symmetric repeated game where the players do not learn anything.

Given a positive integer $n$, the $n$-stage game is the game whose payoff is the expected average  payoff during the first $n$ stages. Under mild assumptions, it has a value, denoted $v_n$. One strand of the literature studies the convergence of the sequence of $n$-stage values, $(v_n)_{n \geq 1}$, as $n$ goes to infinity.

%Given $\lambda \in (0,1]$, the discounted game is the game whose payoff is the discounted sum of the stage payoffs with weight $\lambda$.
%
%convergence the value of discounted game, where the payoff is , when the players become more and more patient.
%
%In order to study the infinitely repeated game a first approach, called \emph{asymptotic} approach, is to study the convergence of the sequence of $n$-stage values, $(v_n)_{n \geq 1}$. Usually one also study discounted games, where the payoff is the discounted sum of the stage payoffs and the convergence of the discounted values when the players become more patient.

The convergence of the sequence of $n$-stage values is related to the behavior of the infinitely repeated game. If the sequence of $n$-stage values converges to some real number $v^*$, one may also consider the existence of a strategy that yields a payoff close to $v^*$ in every sufficiently long game. Let $v$ be a real number. A strategy of player $1$ \emph{guarantees} $v$ if for every $ \eta> 0$ the expected average payoff in the $n$-stage game is greater than $v-\eta$ for every sufficiently large $n$ and every strategy of player $2$. Symmetrically, a strategy of player $2$ \emph{guarantees} $v$ if for every $ \eta> 0$ the expected average  payoff in the $n$-stage game is smaller than $v+\eta$ for every sufficiently large $n$ and every strategy of player $1$. If for every $\epsilon>0$, player $1$ has a strategy that guarantees $v^*-\epsilon$ and player $2$ has a strategy that guarantees $v^*+\epsilon$, then $v^*$ is the \emph{uniform value}. Informally, the players do not need to know the length of the game to play well, provided the game is long enough.

At each stage, the players are allowed to chose their action randomly. If each player can guarantee $v^*$ while choosing at every stage one action instead of a probability over his actions, we say that the game has \emph{a uniform value in pure strategies}.\\

%Given $\lambda \in (0,1]$, the discounted game is the game whose payoff is the discounted sum of the stage payoffs with weight $\lambda$. Whenever the uniform value exists, then the sequence of  $n$-stage value converges to $v^*$ and the disvoun

Blackwell \cite{Blackwell_62} proved that an MDP with a finite set of states and a finite set of actions has a uniform value $v^*$  and the decision maker has a pure strategy that guarantees $v^*$. Moreover, at every stage, the optimal action depends only on the current state. Dynkin and Ju{\v{s}}kevi{\v{c}} \cite{Dynkin_79}, and Renault \cite{Renault_2011} described sufficient conditions for the
existence of the uniform value when the set of states is compact, but in this more general setup there may not exist a strategy which guarantees the uniform value.

Rosenberg, Solan, and Vieille \cite{Rosenberg_2002} proved that POMDPs with a finite
set of states, a finite set of actions, and a finite set of signals
have a uniform value. Moreover, for any $\varepsilon>0$, a
strategy that guarantees $v^*-\varepsilon$ also yields a payoff close to $v^*$ for other
criteria of evaluation like the discounted evaluation. The existence
of the uniform value was extended by Renault \cite{Renault_2011} to
any space of actions and signals, provided that at each stage only a finite number of signals can be realized.
% the
%decision maker chooses a random distribution with finite support
%over his actions and .

In the framework of two-player games, Mertens and Neyman \cite{Mertens_81} showed the existence of the uniform value in stochastic games with a finite set of states and finite sets of actions. Their proof relies on an algebraic argument using the finiteness assumptions. Renault \cite{Renault_2012a} proved the existence of the uniform value for a two-player game where one player controls the transition and the set of states is a compact subset of a normed vector space.

There is an extensive literature about repeated games in which the players are not perfectly informed about the state or the actions played. Rosenberg, Solan and Vieille \cite{Rosenberg_2009} showed, for example, the existence of the uniform value in some class of games where the players observe neither the states nor the actions played by the other players. Another particular class that is closely related to the model we consider is repeated games with symmetric signals. At each stage, the players observe the past actions played and a public signal. Kohlberg and Zamir \cite{Kohlberg_74b} and Forges \cite{Forges_82}
proved the existence of the uniform value, when the state is fixed once and for all at the outset of the game. Neyman and Sorin \cite{Neyman_98} extended this result to the non
zero-sum case, and Geitner \cite{Geitner_2002} to an intermediate model where the game is half-stochastic and half-repeated.  In these four papers, the unknown information concerns a parameter which does not change during the game. We will study models where this parameter can change during the game. \\

In  general, the uniform value does not exist if the players have different information. For example, repeated games with incomplete information on both sides do not have a uniform value \cite{Aumann_95}. Nevertheless, Mertens and Zamir \cite{Mertens_71}
\cite{Mertens_80} showed that the sequence of $n$-stage values converges. Rosenberg, Solan, and Vieille \cite{Rosenberg_2003} and
Coulomb \cite{Coulomb_2003} showed the existence of two quantities,
called the $\max\min$ and the $\min\max$, when each player observes
the state and his own actions but has imperfect monitoring of the
actions of the other player. Moreover the $\max\min$ where player $2$ chooses his strategy knowing the strategy
of player $1$, only depends on the information of player $1$.

More surprisingly, the sequence $(v_n)_{n\geq1}$ may not converge even with symmetric information. Vigeral \cite{Vigeral_2013} provided an example of a stochastic game with a finite set of states and compact sets of actions where the sequence of $n$-stage values does not converge. Ziliotto \cite{Ziliotto_2013} provided an example of a repeated game with a finite set of states, finite sets of actions, and a finite set of public signals where a similar phenomenon occurs. In each case, the game has no uniform value.\\

In this paper, we are interested in two-player zero-sum stochastic games where the
transition is commutative. In such games, given a sequence of decisions, the order is irrelevant to
determine the resulting state: playing an action profile $a_1$ followed by an
action profile $a_2$ leads to the same distribution over states as
playing first the action profile $a_2$ and then $a_1$.

%A Markov chain
%for example can be interpreted as a stochastic game where each
%player has a unique action and the commutativity assumption is
%automatically fulfilled.
%
%The exploitation of a mineral resource such as oil or gold is an
%example of an economic problem fitting this assumption:  it is
%enough to remember how much of the resource has been exploited in
%the past to define the remaining quantity. Another example is a
%competition between firms, which have to sell some stocks. If we
%consider the vector of stocks of all firms as the state and the
%quantities sold by each firm at each stage as the actions, the new
%state depends on the past decisions but not on their order.
%Nevertheless the stage payoff depends on the state and on the
%actions, thus two profiles of actions may lead to different payoffs
%depending on the order they are played.

In game theory, several models satisfy this assumption. For
example, Aumann and Maschler \cite{Aumann_95} studied repeated games with incomplete information on one side. One can introduce an auxiliary
stochastic game where the new state space is the set
of beliefs of the uninformed player and the sets of actions are
the mixed actions of the original game. This game is commutative as we will show in Example \ref{aumann}.
We will also show in Proposition \ref{absorbant} that absorbing games  (see Kohlberg \cite{Kohlberg_74a}) can be reformulated as
commutative stochastic games.

In Theorem \ref{theo1} we prove that whenever a commutative MDP with a deterministic transition has a uniform value in pure strategies, the decision maker has a strategy that guarantees the value.
Example \ref{escalier} shows that to guarantee the value, the decision maker may need to choose his actions randomly. Under topological assumptions similar to Renault
\cite{Renault_2011}, we show that the conclusion can be strengthened
to the existence of a strategy without randomization that guarantees the value.
%%We will study on MDPs with a finite set of actions and a compact state space where the transition is commutative and deterministic. In theorem $\ref{theo1}$, we prove that whenever such a MDP has a uniform value in pure strategies, the decision maker has a strategy which guarantees exactly the value. Without tx�he commutation assumption, this result is false since the decision maker may have to make small irreversible mistakes. Moreover under topological assumptions similar to Renault \cite{Renault_2011} which ensures the existence of the uniform value in pure strategies, we can build this strategy without randomization.
By  a standard argument, we deduce the existence of a strategy that guarantees the value in commutative POMDPs where the decision maker has no information on the state.

%and the transition is commutative, can  be reduced to study an auxiliary MDP with a commutative and deterministic transition. Therefore, we can
%apply Theorem $\ref{theo1}$ to the auxiliary MDP and deduce the
%existence of strategies which guarantee the uniform value.
%The usual approach is to make some assumptions on the signalling structure in order to state a stochastic game on an auxiliary state space as for MDP with Partial Observations and to make some assumptions on the transitions in order to study this auxiliary stochastic game. Indeed this game is not defined on a finite state space so we do not know in general if there exist a uniform value. Then we deduce some result on the original stochastic game with signals.
%For two-player zero-sum stochastic games, we also have a notion of
%uniform value. The

In Theorem \ref{theo2} we prove that a two-player zero-sum stochastic game in which the set of states is a compact subset of $\R^m$, the sets of actions are finite, and each
transition is a deterministic function that is $1$-Lipschitz for the norm
$\|.\|_1$, has a uniform value. We deduce the existence of the uniform value in commutative \emph{state-blind} repeated games where at each stage  the players learn the past actions played but not the state. In this case, we can define an auxiliary
stochastic game on a compact set of states, which satisfies
the assumptions of Theorem \ref{theo2}. Therefore, this auxiliary
game has a uniform value and we deduce the existence of the uniform value in the original
state-blind repeated game.\\

The paper is organized as follows. In Section $2$, we introduce the formal definition of commutativity,
the model of stochastic games, and the model of state-blind
repeated games. In Section $3$, we state the results. Section $4$ is
dedicated to several results on Markov Decision Processes. We first provide an example of a commutative deterministic MDP
with a uniform value in pure strategies but no pure $0$-optimal
strategies. Then we prove Theorem $\ref{theo1}.$ In Section $5$, we focus on the results in
the framework of stochastic games and the proof of Theorem $\ref{theo2}$. We first show
that another widely studied class of games, called absorbing game,
can be reformulated into the class of commutative games. Then we prove Theorem
$\ref{theo2}$ and deduce the existence of the uniform value in commutative state-blind repeated games. Finally, we provide some extensions of Theorem $\ref{theo2}$. Especially, we show the existence of a uniform equilibrium in two-player non zero-sum state-blind commutative repeated games and $m$-player state-blind product-state commutative repeated games.

\section{The model}

%\subsection{Stochastic games.}

When $X$ is a non-empty set, we denote by $\Delta_f(X)$ the set of
probabilities on $X$ with finite support. When $X$ is finite, we
denote the set of probabilities on $X$ by $\Delta(X)$ and by $\sharp
X$ the cardinality of $X$. We will consider two types of games:
stochastic games on a compact metric set $X$ of states, denoted by
$\Gamma=(X,I,J,q,g)$, and state-blind repeated games on a finite
set $K$ of states\footnote{We use $X$ to denote a general set of states and $K$ to denote a finite set of states}, denoted by $\Gamma^{sb}=(K,I,J,q,g)$. The sets of
actions will always be finite. Finite sets will be given the discrete topology. We first define
stochastic games and the notion of uniform value. We will then describe state-blind repeated games.
%When the second player has no influence on the payoff and on the transition, we will forget the set $J$ and called the problem respectively a MDP and a MDP in the dark.

\subsection{Commutative stochastic games}

A two-player zero-sum stochastic game $\Gamma=(X,I,J,q,g)$ is given by a non-empty set of states $X$, two finite, non-empty sets of actions $I$ and $J$, a reward function $g: X \times I \times J \rightarrow [0,1]$ and a transition function $q : X \times I \times J \rightarrow \Delta_f(X) $.\\

Given an initial probability distribution $z_1\in \Delta_f(X)$, the
game $\Gamma(z_1)$ is played as follows. An initial state $x_1$ is
drawn according to $z_1$ and announced to the players. At each stage
$t\geq 1$, player $1$ and player $2$ choose simultaneously
actions, $i_t \in I$  and $j_t \in J$. Player $2$ pays to player $1$
the amount $g(x_t,i_t,j_t)$ and a new state $x_{t+1}$ is drawn
according to the probability distribution $q(x_t,i_t,j_t)$. Then,
both players observe the action pair $(i_t,j_t)$ and the state
$x_{t+1}$. The game proceeds to stage $t+1$. When the initial
distribution is a Dirac mass at $x_1\in X$, denoted by $\delta_{x_1}$, we
denote by $\Gamma(x_1)$ the game $\Gamma(\delta_{x_1})$.\\

If for every initial state and every action pair, the image of $q$ is a Dirac measure, $q$ is said to be \emph{deterministic}.\\

Note that we assume that the transition maps to the set of
probabilities with finite support on $X$: given a stage, a state and
an action pair, there exists a finite number of states possible at the next stage. We equip $X$ with any $\sigma$-algebra $\mathcal{X}$ that includes all countable sets. When $(X,d)$ is a metric space, the Borel $\sigma$-algebra suffices.\\

%Moreover the transition is commutative if the state is invariant under permutation of past couple of actions.
For all $i\in I$ and $j\in J$ we extend $q(\cdot,i,j)$ and $g(\cdot,i,j)$
linearly to $\Delta_f(X)$ by
\begin{align*}
\forall z\in \Delta_f(X), \   \widetilde{q}(z,i,j)=\sum_{x\in X} z(x)q(x,i,j) \text{  and  }  \widetilde{g}(z,i,j)=\sum_{x\in X} z(x)g(x,i,j).
\end{align*}

%If all images are Dirac measures, the transition mapping $q$ is deterministic and we call the unique sequence of states with positive probability the play.

\begin{definition}
The transition $q$ is commutative on $X$ if for all $x \in X$, for all $i,i' \in I$ and for all $j,j'\in J$,
\[
\widetilde{q}(q(x,i,j),i',j')=\widetilde{q}(q(x,i',j'),i,j).
\]
%where $q(.,i,j)$ and $q(.,i',j')$ are linearly extended to $\Delta_f(X)$.
\end{definition}

\noindent That is, the distribution over the state after two
stages is equal whether action pair $(i,j)$ is played before
$(i',j')$ or whether $(i,j)$ is played after $(i',j')$. Note that if
the transition $q$ is not deterministic, $\widetilde{q}(q(x,i',j'),i,j)$
is the law of a random variable $x''$ computed in two steps: $x'$ is
randomly chosen with law $q(x,i',j')$, then $x''$ is randomly chosen
with law $q(x',i,j)$; specifically $(i,j)$ is played at the second step independently of the outcome of $x'.$

\begin{remark}
\rm If the transition $q$ does not depend on the actions, then the state process is a Markov chain and the
commutativity assumption is automatically fulfilled.
\end{remark}

% We give two examples of commutative games.

\begin{example}\label{exemple1}
\rm Let $X$ be the set of complex numbers of modulus $1$ and $\alpha:I \times J \rightarrow \Delta_f([0,2\pi])$. Let $q$ be defined by
\[
\forall x \in X,\forall a\in I ,\forall b\in J, \
q(x,a,b)=\sum_{\rho \in[0,2\pi]}\alpha(a,b)(\rho)\delta_{xe^{i \rho}}.
\]

If the state is $x$ and the action pair $(a,b)$ is played,
then the new state is $x'=xe^{i \rho}$ with probability
$f(a,b)(\rho)$. This transition is commutative by the commutativity of multiplication of complex numbers.
\end{example}

The next example originates in the theory of repeated games with incomplete information on one side (Aumann and Maschler \cite{Aumann_95}).
%The auxiliary stochastic game associated to the repeated game is commutative.

% The revelation of information from player $1$ to player $2$ is commutative,

\begin{example}\label{aumann}
\rm A repeated game with incomplete information on one side,
$\Gamma$, is defined by a finite family of matrices $(G^k)_{k\in
K}$, two finite sets of actions $I$ and $J$, and an initial
probability $p_1$. At the outset of the game, a matrix $G^k$ is randomly chosen
with law $p_1$ and told to player $1$ whereas player $2$ only knows
$p_1$. Then, the matrix game $G^k$ is repeated over and over. The
players observe the actions played but not the payoff.

One way to study $\Gamma$ is to introduce a stochastic game on the posterior beliefs of player $2$ about the
state.  Knowing the strategy played by player $1$, player $2$ updates his
posterior belief depending on the actions observed. Let $\Psi=(X,A,B,\widetilde{q},\widetilde{g})$ be a
stochastic game where $X=\Delta(K)$, $A=\Delta(I)^K$ and
$B=\Delta(J)$, the payoff function is given by
\[
 \widetilde{g}(p,a,b)=\sum_{k\in K,i\in I, j\in J}p^k a^k(i)b(j)G^k(i,j),
\]
 and the transition by
\[
\widetilde{q}(p,a,b)=\sum_{k \in K,i \in I} a^k(i) \delta_{\hat{p}(a|i)},
\]
where $a(i)=\sum_{k\in K} p^k a^k(i)$ and $\hat{p}(a|i)=\left( \frac{p^k a^k(i)}{a(i)}\right)_{k\in K} \in \Delta(K)$.
Knowing the mixed action chosen by player $1$ in each state, $a$, and having a prior belief, $p$, player $2$ observes action $i$ with probability $a(i)$ and updates his beliefs by Bayes rule to $\hat{p}(a|i)$. This induces the auxiliary transition $\widetilde{q}.$ The payoff $\widetilde{g}$ is the expectation of the payoff under the probability generated by player 2's belief and player 1's mixed action.

We now check that the auxiliary  stochastic game is commutative. Note that the
second player does not influence the transition so we can ignore him. Let $a$ and $a'$ be two actions of player $1$ and $p$ be a belief of player $2$. If player $1$ plays first $a$ and player $2$ observes action $i$, then
player $2$'s belief $p_2(\cdot|i)$ is given by
\[\forall k\in K, \ p_2(k|i)=\frac{p^k a^k (i)}{\sum_{k\in K} p^k
a^k(i)}.
\]
If now player $1$ plays $a'$ and player $2$ observes $i'$, then player $2$'s belief $p_3(\cdot|i,i')$ is given by
\[
\forall k\in K, \ p_3(k|i,i')=\frac{p_2(k|i)a'^k(i')}{\sum_{k\in K} p_2(k|i)a'^k(i')}=\frac{p^k a'^k(i')a^k(i)}{\sum_{k\in K} p^k a'^k(i')a^k(i)}.
\]
The probability that the action pair $(i,i')$ is observed is $p^k
a^k(i) a'^k(i')$. Since the belief $p_3$ and the probability to observe each pair $(i,i')$ are symmetric in $(a,i)$ and $(a',i')$, the transition is commutative.
\end{example}

\begin{remark}
\rm Both previous examples are commutative but the transition is not
deterministic.
\end{remark}

\begin{remark}
\rm Commutativity of the transitions implies that if we consider an initial state $x$ and a finite
sequence of actions $(i_1,j_1,....,i_n,j_n)$, then the law of the state
at stage $n+1$ does not depend on the order in which the action pairs $(i_t,j_t)$, $t=1,...,n,$ are played. We can thus represent a finite sequence of actions by a vector in $\N^{I\times J}$ counting
how many times each action pair is played. Other models in which the transition
along a sequence of actions is only a function of a parameter in a
smaller set have been studied in the literature. For example, a transition is state independent (SIT) if
it does not depend on the state. The law of the state at stage $n$
is characterized only by the last action pair played. The law then
depends on the order in which actions are played.
Thuijsman \cite{Thuijsman_92} proved the existence
of stationary optimal strategies in this framework.
\end{remark}

\subsection{Uniform value}\label{stoc}

%Note that given an initial probability $z\in \Delta_f(X)$ with finite support, q since there are
%finitely many couples of actions, the set of states, which can be
%reached with positive probability, is countable. Therefore, the
%actual choice of algebra is not important.

At stage $t$, the space of past \emph{histories} is
$H_{t}=(X \times I \times J)^{t-1}\times X.$ Set
$H_{\infty}=(X\times I \times J)^{\infty}$  to be the space of infinite
\emph{plays}. For every $t\geq 1$, we consider the product topology on $H_t$, $t\geq 1$ and also on $H_\infty$.\\

%Without additional assumption on $X$, there exists an
%infinite number of histories, possibly uncountable, in $H_t$ for
%every $t\geq 1$ and in $H_\infty$.\\

A \emph{(behavioral) strategy} for player $1$ is a sequence $(\sigma_t)_{t
\geq 1}$ of functions $\sigma_t: H_t \rightarrow \Delta(I)$. A
\emph{(behavorial) strategy} for player $2$ is a sequence $\tau=(\tau_t)_{t
\geq 1}$ of functions $\tau_t: H_t \rightarrow \Delta(J)$. We denote
by $\Sigma$ and $\T$, the player's respective sets of strategies.

Note that we did not make any measurability assumption on the strategies. Given $x_1\in X$, the set of histories at stage $t$ from state $x_1$ is finite since the image of the transition $q$ is contained in the set of probabilities over $X$ with finite support and the sets of actions are finite. It follows that any triplet $(z_1,\sigma,\tau)$ defines a probability over $H_t$ without an additional measurability condition. This sequence of probabilities  can be extended to a unique probability denoted $\PP_{z_1,\sigma,\tau}$ over the set $H_\infty$ with the infinite product topology. We denote by $\E_{z_1,\sigma,\tau}$ the expectation with respect to the probability $\PP_{z_1,\sigma,\tau}$.

%Since there exists only a finite number of actions and the
%following, we can restrict without loss of generality to a countable
%subset of $H_\infty$. We consider the product topology on $H_t$ and
%the Borel $\sigma$-algebra associated, which contains all countable
%sets and complements of countable sets. Each $h_t \in H_t$ can be
%identified with a cylinder set of $H_{\infty}$. We denote by
%$\mathcal{H}_t$ the algebra induced by $H_t$ over $H_{\infty}$ and
%by $\mathcal{H}_{\infty}$ the $\sigma$-algebra spanned by all finite
%cylinders. We do not makeThe set of actions is finite and $q$ has
%values in laws with finite support, so we do not need assumption of
%measurability.
%
%Let $(\sigma,\tau)$ be a strategy profile and $z$ be an initial
%probability with finite support on $X$,
% induce a probability on each finite cylinder, which can be
%extended as a unique probability distribution $\P_{z,\sigma,\tau}$
%over the set of infinite histories
%$(H_{\infty},\mathcal{H}_{\infty})$. The set of histories with
%positive probability is countable.

If for every $t\geq 1$ and every history $h\in H_t$ the image of $\sigma_t(h_t)$ is a Dirac measure, the
strategy is said to be \emph{pure}. If the initial distribution is a Dirac measure, the transition is
deterministic and both players use pure strategies, then $\PP_{z_1,\sigma,\tau}$ is a Dirac measure. The strategies induce a unique play.

The game we described is a game with perfect recall, so that by Kuhn's theorem \cite{Kuhn_53} every behavior strategy is equivalent to a probability over pure strategies, called \emph{mixed strategy}, and vice versa.\\
%one history with a positive probability. We call it the play and it
%can be uniquely defined either by the sequence of states visited or
%by the sequence of actions played.
%
%Another way to define a strategy is by the introduction of mixed
%strategies, that is, probability distribution over the pure
%strategy. By Kuhn's lemma, the sets of probabilities generated on
%histories by behavioral strategies and mixed strategies are the
%same.\\

We are going to focus on two types of evaluations, the $n$-stage expected payoff
and the expected average payoff between two stages $m$ and $n$. For
each positive integer $n$, the expected average payoff for player $1$ up to
stage $n$, induced by the strategy pair $(\sigma,\tau)$ and the
initial distribution $z_1$, is given by
\[
\gamma_n(z_1,\sigma,\tau)=\E_{z_1,\sigma,\tau}\left(\frac{1}{n} \sum_{t=1}^n g(x_t,i_t,j_t)\right).
\]
\noindent The expected average payoff between two stages $1 \leq m \leq n$ is given by
\[
\gamma_{m,n}(z_1,\sigma,\tau)=\E_{z_1,\sigma,\tau}\left(\frac{1}{n-m+1} \sum_{t=m}^{n} g(x_t,i_t,j_t)\right).
\]

\noindent To study the infinitely repeated game $\Gamma(z_1)$, we focus on the notion of uniform value and on the notion of $\varepsilon$-optimal strategies.
%We use a more general definition than usual by restricting the available strategies for the players.

\begin{definition}\label{stoc_uni}
Let $v$ be a real number.
\begin{itemize}
\item Player $1$ can \emph{guarantee} $v$ in $\Gamma(z_1)$ if for all $\varepsilon >0$ there exists a strategy $\sigma^*\in \Sigma$ of player $1$ such that
\[
\liminf_n \inf_{\tau \in \T} \gamma_n(z_1,\sigma^*,\tau) \geq v-\varepsilon.
\]
We say that such a strategy $\sigma^*$ \emph{guarantees} $v-\varepsilon$ in $\Gamma(z_1)$.
\item Player $2$ can \emph{guarantee} $v$ in $\Gamma(z_1)$ if for all $\varepsilon >0$ there exists a strategy $\tau^* \in \T$ of player $2$ such that
\[
\limsup_n \sup_{\sigma \in \Sigma} \gamma_n(z_1,\sigma,\tau^*) \leq v+\varepsilon.
\]
We say that such a strategy $\tau^*$ \emph{guarantees} $v+\varepsilon$ in $\Gamma(z_1)$.

\item If both players can guarantee $v$, then $v$ is called the \emph{uniform value} of the game $\Gamma(z_1)$ and denoted by $v^*(z_1)$. A strategy $\sigma$ (resp. $\tau$) that guarantees $v^*(z_1)-\varepsilon$ (resp. $v^*(z_1)+\epsilon$) with $\varepsilon \geq 0$ is called \emph{$\varepsilon$-optimal}.
\end{itemize}

\end{definition}

\begin{remark}
\rm For each $n\geq 1$ the triplet
$(\Sigma,\T,\gamma_n(z_1,.,.))$ defines a game in strategic form. This game has a
value, denoted by $v_n(z_1)$. If the game
$\Gamma(z)$ has a uniform value $v^*(z_1)$, then the sequence $(v_n(z_1))_{n\geq 1}$
converges to $v^*(z_1)$.
\end{remark}

\begin{remark}
\rm Let us make several remarks on another way to evaluate the infinite stream of payoffs. Let $\lambda \in(0,1]$. The expected $\lambda$-discounted payoff for
player $1$, induced by a strategy pair $(\sigma,\tau)$ and the
initial distribution $z_1$, is given by
\[
\gamma_\lambda(z_1,\sigma,\tau)=\E_{z_1,\sigma,\tau}\left(\lambda
\sum_{t=1}^\infty (1-\lambda)^{(t-1)} g(x_t,i_t,j_t)\right).
\]
For each $\lambda \in (0,1]$, the triplet
$(\Sigma,\T,\gamma_\lambda(z_1,.,.))$ also defines a game $\Gamma_\lambda(z_1)$ in strategic form. The sets of strategies are compact for the product topology, and the payoff function $\gamma_\lambda(z_1,\sigma,\tau)$ is continuous. Using Kuhn's theorem, the payoff is also concave-like, convex-like and it  follows therefore from Fan's minimax theorem (see \cite{Fan_53}) that the game $\Gamma_\lambda(z_1)$ has a value, denoted $v_\lambda(z_1)$. Note that there may not exist an optimal measurable strategy which depends only on the current state (Levy \cite{Levy_2012}).

Some authors focus on the existence of $v(z_1)$ such that
\[
\lim_{n \rightarrow \infty} v_n(z_1)=\lim_{\lambda \rightarrow 0}
v_\lambda(z_1)=v(z_1).\]
When the uniform value exists, this equality is immediately true with $v(z_1)=v^*(z_1)$
since the discounted payoff can be written as a convex combination
of expected average payoffs.

\end{remark}
%
%\begin{remark}
%
%\end{remark}

%\noindent When $\Sigma'=\Sigma$ and $\T'=\T$, it is the classic definition of the uniform value.

%Whenever the uniform value exists,  A strategy $\tau$ of player $2$ is
%\emph{$\varepsilon$-optimal} if it guarantees $v^*(z)+\varepsilon$. The
%uniform value exists if and only if there exists $v$ such that
%player $1$ and player $2$ have strategies,
%which guarantee respectively $v-\varepsilon$ and $v+\varepsilon$ for each $\varepsilon>0$.

%We denote by $\max \min$ the supremum of the values that player $1$ can guarantee and $\min \max$ the infimum of the values that player $2$ can guarantee. It is easy to see that the game has a uniform value if both are equal.

\subsection{The model of repeated games with state-blind players}
A \emph{state-blind} repeated game $\Gamma^{sb}=(K,I,J,q,g)$ is defined by the same objects as a stochastic game. The definition of commutativity is the same. The main difference is the information that the players have, which affects their sets of strategies. We assume that at each stage, the players observe the actions played but not the state. We will restrict the discussion to a finite state space $K$.

Given an initial probability $p_1\in \Delta(K)$, the game
$\Gamma^{sb}(p_1)$ is played as follows. An initial state $k_1$ is
drawn according to $p_1$ without being announced to the players. At
each stage $t\geq 1$, player $1$ and player $2$ choose
simultaneously an action, $i_t \in I$  and $j_t \in J$. Player $1$
receives the (unobserved) payoff $g(k_t,i_t,j_t)$, player $2$
receives the (unobserved) payoff $-g(k_t,i_t,j_t)$, and a new state
$k_{t+1}$ is drawn according to the probability distribution
$q(k_t,i_t,j_t)$. Both players then observe only the action pair $(i_t,j_t)$ and the game proceeds to stage $t+1$.

Since the states are not observed, the space of public histories of length $t$ is
$H^{sb}_{t}=(I \times J)^{t-1}$. A
strategy of player $1$ in $\Gamma^{sb}$ is a sequence $(\sigma_t)_{t
\geq 1}$ of functions $\sigma_t: H^{sb}_t \rightarrow \Delta(I)$, and a
strategy of player $2$ is a sequence $\tau=(\tau_t)_{t \geq 1}$ of
functions $\tau_t: H^{sb}_t \rightarrow \Delta(J)$. We denote by
$\Sigma^{sb}$ and $\T^{sb}$ the players respective sets of strategies. An
initial distribution $p_1$ and a pair of strategies
$(\sigma,\tau)\in \Sigma^{sb} \times \T^{sb}$ induce a unique
probability over the infinite plays $H_{\infty}$. For every pair of strategies $(\sigma,\tau)$ and initial probability $p_1$ the payoff is defined as in Section \ref{stoc}. Similarly, the notion of uniform value is defined as in Definition \ref{stoc_uni} by restricting the players to play strategies in $\Sigma^{sb}$ and $\T^{sb}$.

\begin{definition}
Let $v$ be a real number.
\begin{itemize}
\item Player $1$ can \emph{guarantee} $v$ in $\Gamma^{sb}(p_1)$ if for all $\varepsilon >0$ there exists a strategy $\sigma^*\in \Sigma^{sb}$ of player $1$ such that
\[
\liminf_n \inf_{\tau \in \T^{sb}} \gamma_n(p_1,\sigma^*,\tau) \geq v-\varepsilon.
\]
We say that such a strategy $\sigma^*$ \emph{guarantees} $v-\varepsilon$ in $\Gamma^{sb}(p_1)$.
\item Player $2$ can \emph{guarantee} $v$ in $\Gamma^{sb}(p_1)$ if for all $\varepsilon >0$ there exists a strategy $\tau^* \in \T^{sb}$ of player $2$ such that
\[
\limsup_n \sup_{\sigma \in \Sigma^{sb}} \gamma_n(p_1,\sigma,\tau^*) \leq v+\varepsilon.
\]
We say that such a strategy $\tau^*$ \emph{guarantees} $v+\varepsilon$ in $\Gamma^{sb}(p_1)$.

\item If both players can guarantee $v$, then $v$ is called the \emph{uniform value} of the game $\Gamma^{sb}(p_1)$ and denoted by $v^{sb}(p_1)$.
\end{itemize}
\end{definition}

%Formally a state-blind general repeated games is a general repeated games $\Gamma^sb=(K,I,J,C,D,q,r)$ where $C$ and $D$ are the signals the player receives and are reduced to singleton. Furthermore $\Sigma^{sb}$ and $\T^{sb}$ are then respectively subsets of $\Sigma$ and $\T$, thus we can also consider this uniform value as a uniform value in the original problem where the sets of strategies are restricted. Note that here is no relation, a priori, between $v^{sb}(p_1)$ and $v(p_1)$ since both sets of strategies are restricted.   and we can define the notion of uniform value.
%
\begin{remark}
\rm The sets $\Sigma^{sb}$ and $\T^{sb}$ can be seen as subsets of $\Sigma$ and $\T$ respectively. There is no relation between $v^{sb}(p_1)$ and $v^*(p_1)$, since both players have restricted sets of strategies.
\end{remark}
%and we can define the notion of uniform value.

%\subsection{Commutative stochastic games}

%%%%%%%%%%%%%%%%%%%%
%  TOPOLOGIE ???
%%

%Furthermore, the set of state on this history will be called the path of the strategy.

%x%
\section{Results.}

In this section we present the main results of the paper. Section
\ref{results_mdp} concerns MDPs and Section \ref{results_stoc}
concerns stochastic games.

\subsection{Existence of $0$-optimal strategies in Commutative deterministic Markov Decision Processes.}\label{results_mdp}

An MDP is a stochastic game, $\Gamma=(X,I,q,g)$, with a single player, that is, the set $J$ is a singleton. Our first main result states that if an MDP with deterministic and commutative transitions has a uniform value and if the decision maker has pure $\epsilon$-optimal strategies, then he also has a (not necessarily pure) $0$-optimal strategy. We also provide sufficient topological conditions for the existence of a pure $0$-optimal strategy.

\begin{theorem}\label{theo1}
Let $\Gamma=(X,I,q,g)$ be an MDP such that $I$ is finite and $q$ is deterministic and commutative.
\begin{enumerate}
\item If for all $z_1\in \Delta_f(X)$, $\Gamma(z_1)$ has a uniform value in pure strategies, then for all $z_1\in \Delta_f(X)$ there exists a $0$-optimal strategy.
\item If $X$ is a precompact metric space, $q(\cdot,i)$ is $1$-Lipschitz  for every $i\in I$, and $g(\cdot,i)$ is uniformly continuous for every $i\in I$, then for all $z_1\in \Delta_f(X)$ the game $\Gamma(z_1)$ has a uniform value and there exists a $0$-optimal pure strategy.
\end{enumerate}
\end{theorem}

\begin{remark}
\rm In an MDP with a deterministic transition, a play is uniquely determined by the initial state and a sequence of actions. Thus, in the framework of deterministic MDPs we will always identify the set of pure strategies with the set of sequences of actions.
\end{remark}

The first part of Theorem \ref{theo1} is sharp in the sense that a
commutative deterministic MDP with a uniform value in pure
strategies may have no $0$-optimal pure strategy. An example is
described at the beginning of Section $4$.

The topological assumptions of the second part of Theorem \ref{theo1} were first introduced
by Renault \cite{Renault_2011} and imply the existence of the
uniform value in pure strategies; by the first part of the theorem they also imply the existence of a $0$-optimal strategy. Under these topological
assumptions, we prove the stronger result of the existence of a
$0$-optimal pure strategy.

%The decision maker can guarantee the
%payoff exactly without randomizing.

Let us now discuss the topological assumptions made in Theorem \ref{theo1}. First, if the
payoff function $g$ is only continuous or the state space  is not precompact,
then the uniform value may fail to exist as shown in the following
example.

\begin{example}
\rm Consider a Markov Decision Process $(X,I,q,g)$ where there is only one
action, $|I|=1$. The set of states is the set of integers, $X=\N$, and the
transition is given by $q(n)=n+1, \ \forall n\in \N$. Note that $q$ is commutative and deterministic. Let $r=(r_n)_{n \in \N}$ be a sequence of numbers in $[0,1]$ such that the sequence of average payoffs does
not converge. The payoff function is defined by $g(n)=r_n, \forall n\in \N$.

We consider the following metric on $\N:$ for all $n,m \in \N$,
$d(n,m)=\1_{n\neq m}$. Then $(\N,d)$ is not precompact, the transition $q$ is $1$-Lipschitz, and the
function $g$ is uniformly continuous. The choice of $r$ implies that the MDP $\Gamma=(X,I,q,g)$ has no uniform value.

Consider now the following metric on $\N:$ for all $n,m \in \N$,
$d'(n,m)=|\frac{1}{n+1}-\frac{1}{m+1}|$. Then $(\N,d')$ is a
precompact metric space, the transition is $1$-Lipschitz and the function $g$ is continuous. As before, the MDP $\Gamma=(X,I,q,g)$ has no uniform value. A simple computation shows that the function $g$ is not uniformly
continuous on $(\N,d')$. Take now $g$ a uniformly continuous function, then $(g(n))_{n \in \N}$ is a Cauchy sequence in a complete space, thus converges. It follows that the sequence of Ces\`aro averages also converges to the same limit and the game has a uniform value.
\end{example}

The assumption that $q$ is $1$-Lipschitz may seem strong but turns out to be necessary in the proof of Renault \cite{Renault_2011}. The reason is as follows. When computing the uniform value, one considers
infinite histories. When $q$ is $1$-Lipschitz, given two states $x$
and $x'$ and an infinite sequence of actions $(i_1,...,i_t,...)$, the state at stage $t$ on the play from $x$ and the state at stage $t$ on the play
from $x'$ are at a distance at most $d(x,x')$. Thus, the payoffs
along both plays stay close at every stage. On the contrary, if $q$
were say $2$-Lipschitz, we only know that the distance between the
state at stage $t$ on the play from $x$ and the state at stage $t$ on the play from $x'$ is
at most  $d(x,x') 2^t$, which gives no uniform bound on the difference between the stage payoffs along the two plays. As shown in Renault \cite{Renault_2011} when $q$ is not $1$-Lipschitz, the value may fail to exist. The counter-example provided by Renault is not commutative and it might be that the additional assumption of commutativity can
help us in relaxing the Lipschitz requirement on $q$. In our proof, we use the fact that $q$
is $1$-Lipschitz at two steps: first in order to apply the result of
Renault \cite{Renault_2011} and then in order to concatenate
strategies. It is still open whether one of these two steps can be done under
the weaker assumption that $q$ is uniformly continuous. 

We list now two open problems: assume that the uniform value exists, $X$ is
precompact, $g$ is uniformly continuous, and $q$ is uniformly
continuous, deterministic, and commutative; does there exist a
$0$-optimal strategy? Does an MDP with $X$ precompact, $g$ uniformly
continuous, and $q$ uniformly continuous, deterministic, and
commutative  always
have a uniform value?\\

  %If $i\in I$ and we denote by $q_i=q(.,i)$ the transition if $i$ is played, then $q_i$ $1$-Lipschitz implies that for all $n\geq 1$, $q_i^n$ is also $1$-Lipschitz

We deduce from Theorem \ref{theo1} the existence of a $0$-optimal
strategy for commutative POMDPs with no information on the state,
called \emph{MDPs in the dark} in the literature. The auxiliary MDP
associated to the POMDP is deterministic and commutative, and thus it
satisfies the assumption of Theorem \ref{theo1}.

\begin{corollary}\label{MDPstateblind}
Let $\Gamma^{sb}=(K,I,q,g)$ be a commutative state-blind POMDP with a finite state space $K$ and a finite set of actions $I$. For all $p_1\in \Delta(K)$, the POMDP $\Gamma^{sb}(p_1)$ has a uniform value and there exists a $0$-optimal pure strategy.
\end{corollary}

We will prove Corollary \ref{MDPstateblind} in the two-player framework.\\
%A MDP in the dark can be written as a decision problem on the space of beliefs with the same set of actions. Since the initial set of states is finite, the space of beliefs is compact and we can apply the theorem.

Rosenberg, Solan, and Vieille \cite{Rosenberg_2002} asked if a $0$-optimal strategy exists in POMDPs. Theorem \ref{theo1} ensures that if the transition is commutative such a strategy exists. The following example, communicated by Hugo Gimbert,  shows that it is not true without the
commutativity assumption. The example also implies that there
exist games that cannot be transformed into a commutative game with
finite sets of actions.

\begin{example}\label{noirRSV}
\rm Consider a state-blind POMDP $\Gamma^{sb}=(X,I,g,q)$ defined as follows. There are four states $X=\{\alpha,\beta,k_0,k_1\}$ and two actions $I=\{T,B\}$. The payoff is $0$ except in state $k_1$ where it is $1$. The states $k_0$ and $k_1$ are absorbing and the transition function $q$ is given on the other states by
\begin{align*}
 q(\alpha,T) & =\frac{1}{2}\delta_{\alpha}+\frac{1}{2}\delta_{\beta} ,\\
 q(\beta,T) & =\delta_{\beta},\\
 q(\alpha,B) & =\delta_{k_0}, \\
 q(\beta,B) & =\delta_{k_1}.\\
\end{align*}
\noindent This POMDP is not commutative: if the initial state is $\alpha$ and the decision maker plays $B$ and then $T$, the state is $k_0$ with probability one, whereas if he plays first $T$ and then $B$, the state is $k_0$ with probability $1/2$ and $k_1$ with probability $1/2$.

Let us check that this game has a uniform value in pure strategies,
but no $0$-optimal strategies. An $\varepsilon$-optimal strategy in
$\Gamma(\alpha)$ is to play the action $T$ until the probability to
be in $\beta$ is more than $1-\varepsilon$, and then to play $B$. This
leads to a payoff of $1-\varepsilon$, so the uniform value in $\alpha$
exists and is equal to $1$. The reader can verify that there is no strategy that guarantees $1$ in $\Gamma(\alpha)$.
%In order to get a good payoff at some
%stage, the decision maker has to play $B$ and reach the state $k_0$
%with positive probability. It is an irreversible mistake in the
%sense that by choosing this action, the expectation of the uniform
%value decreases strictly. From the new distribution, the decision
%maker can no longer guarantee the value $1$ but at most
%$1-\varepsilon$.
\end{example}

\subsection{Existence of the uniform value in commutative deterministic stochastic
games.}\label{results_stoc}

For two-player stochastic games, the commutativity assumption does not
imply the existence of $0$-optimal strategies. Indeed, we will prove in
Proposition \ref{absorbant} that any absorbing game can be
reformulated as a commutative stochastic game. Since there exist
absorbing games with deterministic transitions without $0$-optimal
strategies, for example the Big Match (see Blackwell and Ferguson \cite{Blackwell_68}), there exist deterministic
commutative stochastic games with a uniform value and without $0$-optimal strategies.
In this section, we study the existence of the uniform value in one class of
stochastic games on $\R^m$.

\begin{theorem}\label{theo2}
Let $\Gamma=(X,I,J,q,g)$ be a stochastic game where $X$ is a compact subset of $\R^m$, $I$ and $J$ are finite sets, $q$ is commutative, deterministic and $1$-Lipschitz for the norm $\|.\|_1$, and $g$ is continuous. Then for all $z_1\in \Delta_f(X)$ the stochastic game $\Gamma(z_1)$ has a uniform value.
\end{theorem}

Let us comment on the assumptions of Theorem \ref{theo2}. The state
space is not finite yet the set of actions available to each player
is the same in all states. This requirement is necessary to ensure
that the commutativity property is well defined. Our proof is valid
only if $q$ is $1$-Lipschitz with respect to the norm $\|.\|_1$.
Thus this theorem does not apply to Example \ref{exemple1} on the
circle. The proof can be adapted for polyhedral norms (i.e. such
that the unit ball has a finite number of extreme points), this is
further discussed in Section \ref{extension}. Finally note that the
most restrictive assumptions are on the transition.
%and lead to very specific properties. For example, if $x\in X$ is a fixed point of a commutative deterministic transition $q(.,i^*,j^*)$ with $(i^*,j^*)\in I\times J$, then each state reached in the game $\Gamma(x)$ will still be a fixed point of $q(.,i^*,j^*).$

As shown in the MDP framework, the assumption that $q$ is
$1$-Lipschitz is important for the existence of a uniform value and
is used in the proof at two steps. First, we use it to deduce that
for all $(i,j)\in I\times J$, iterating infinitely often the action pair $(i,j)$ leads to a limit cycle with a finite number of
states. Second, it is used to prove that if a strategy guarantees $w$
from a state $x$ then it guarantees almost $w$ in any game that starts at an initial state in a small neighbourhood of $x$.\\

Given a state-blind repeated game $\Gamma^{sb}=(K,I,J,q,g)$ with a commutative transition $q$, we define the auxiliary stochastic game
$\Psi=(X,I,J,\widetilde{q},\widetilde{g})$ where $X=\Delta(K)$,
$\widetilde{q}$ is the linear extension of $q$, and $\widetilde{g}$ is the
linear extension of $g$.

%A state, in this new game, is the common
%belief of the players over the state in $\Gamma^{sb}(K,I,J,q,g)$.

Since $K$ is finite,  $X$ can be embedded in $\R^{K}$ and the transition $\widetilde{q}$ is deterministic, $1$-Lipschitz for $\|.\|_1$, and commutative. Furthermore, $\widetilde{g}$ is continuous and therefore we can apply Theorem \ref{theo2} to $\Psi$. It follows that for each initial state $p_1 \in X$, $\Psi(p_1)$ has a uniform value. We will check that it is  the uniform value of the state-blind repeated game $\Gamma^{sb}(p_1)$ and deduce the following corollary.

% with full monitoring which fulfills the assumption of the theorem. In $\widetilde{\Gamma}$, given a couple of actions, the transition $\widetilde{q}$ is deterministic so we can focus on strategies which depend only on the actions. This set of strategies is exactly the set of strategies of $\Gamma^{sb}$. Moreover the transition $\widetilde{q}$ is defined such that the state in $\widetilde{\Gamma}$ is the law of the state in $\Gamma^{sb}$. We deduce that the two games have the same finite values and that the existence of a uniform value in $\widetilde{\Gamma}$ implies the existence of the uniform value in the original game.

\begin{corollary}\label{mdp_noir}
Let $\Gamma^{sb}=(K,I,J,q,g)$ be a commutative state-blind repeated game with a finite set of states $K$ and finite sets of actions $I$ and $J$. For all $p_1\in \Delta(K)$, the game $\Gamma^{sb}(p_1)$ has a uniform value.
\end{corollary}

%and for all $n\in \N$, we have $v_n^{sb}(x_1)=\widetilde{v}_n(x_1)$.

\begin{remark}
\rm Corollary \ref{mdp_noir} concerns repeated games where the
players observe past actions but not the state. The more general
model, where the players observe past actions and have a public
signal on the state, leads to the definition of an auxiliary
stochastic game with a random transition. In the
deterministic case given a triplet $(z_1,\sigma,\tau)$, the sequence
of states visited along each infinite play converges to a finite
cycle $\P_{z_1,\sigma,\tau}$-almost surely. This no longer holds if the transition is random.

%In the state-blind case, the commutativity assumption implies that given a pair of strategies, every sequence of states extracted from a play with positive probability is a precompact family.admits a finite number of states
%Indeed, in the definition of commutation, we consider that after
%one stage, the decision maker chooses the same action whatever is
%the intermediate state. Thus the assumption does not give any
%property if the decision maker plays differently depending on the
%signal he has observed after one stage. When the transition is
%deterministic, this problem does not appear, since there is only one
%play with positive probability.
\end{remark}

We now present an example of a commutative state-blind repeated game and its auxilliary deterministic stochastic game.

\begin{example}
\label{exemple:triangle}
\rm Let $K=\Z/m\Z$ and $f$ be a function from  $I \times J$ to $\Delta(K)$. We define the transition $q : K
\times I \times J \rightarrow \Delta(K)$ as follows: given a state $k\in K$, if the players play $(i,j)$ then for all $k'\in K$, the new state is $k+k'$ with probability $f(i,j)(k')$.

If the initial state is drawn with a distribution $p$ and players
play respectively $i$ and $j$, then the new state is given by the
sum of two independent random variables of respective laws $p$ and
$f(i,j)$. The addition of independent random variables is a
commutative and associative operation, therefore $q$ is commutative on
$K$.

For example, let $m=3$, $I=\{T,B\}$, $J=\{L,R\}$ and the function
$f$ be given by
\begin{center}
\begin{tabular}{cc}
 & $\begin{matrix}  L & \hspace{1cm} R \end{matrix} $  \\
$\begin{matrix}
T\\
B
\end{matrix}$ & $\begin{pmatrix}
 \frac{1}{2} \delta_1+\frac{1}{2}\delta_2 & \delta_1 \\
 \delta_1  & \delta_0
\end{pmatrix}$.
\end{tabular}
\end{center}
\noindent If the players play $(T,L)$ then the new state is one of
the two other states with equal probability. If the players play
$(B,R)$, then the state does not change. And otherwise the state
goes from state $k$ to state $k+1$.
%The auxiliary game $\widetilde{\Gamma}=(\Delta(K),I,J,\widetilde{q},\widetilde{r})$ is given by
%\[r(p_0,p_1,p_2)= p_1 r(

The extension of the transition function to the set of probabilities
over $K$ is given by
\begin{center}
\begin{tabular}{rl}
$\widetilde{q}((p^1,p^2,p^3),T,L)$ & $= \left( \frac{p^2+p^3}{2},\frac{p^1+p^3}{2}, \frac{p^1+p^2}{2} \right)$, \\
$\widetilde{q}((p^1,p^2,p^3),B,R)$ & $= \left( p^1, p^2, p^3 \right)$, \\
$\widetilde{q}((p^1,p^2,p^3),B,L)=\widetilde{q}((p^1,p^2,p^3),T,R)$& $= \left( p^3, p^1, p^2 \right)$.
\end{tabular}
\end{center}
\end{example}

%%%%%%%%%%%%%%%%%%%%%%%%%%%%%%%%%%%%%%%%%%%%%%%%%%%%%%%%%%%%%%%%%%%%%%%%%%%%%%%%%%%%%%%%%%%%%%%%%%%%%%%%%%%%%%%
%%%%%%%%%%%%%%%%%%%%%%%%%%%%%%%%%%%%%%%%%%%%%%%%%%%%%%%%%%%%%%%%%%%%%%%%%%%%%%%%%%%%%%%%%%%%%%%%%%%%%%%%%%%%%%%%
%%%%%%%%%%%%%%%%%%%%%%%%%%%%%%%%%%%%%%%%%%%%%%%%%%%%%%%%%%%%%%%%%%%%%%%%%%%%%%%%%%%%%%%%%%%%%%%%%%%%%%%%%%%%%%%

\section{Existence of $0$-optimal strategies in commutative deterministic MDPs.}

%%%%By Kuhn theorem we behavorial and mixed strategies allow us to switch from general strategies to random. We will define our optimal strategy has a probability on the space of pure strategies.

In this section we focus on MDPs and Theorem \ref{theo1}. The section is
divided into four parts. In the first part we provide an example
showing that under the conditions of  Theorem $\ref{theo1}$(1), there need not exist a pure $0$-optimal strategy.

The rest of the section is dedicated to the proof of Theorem
\ref{theo1}. In the second part we show that in a deterministic
commutative MDP, for all $\varepsilon>0$, there exist
$\varepsilon$-optimal pure strategies such that the uniform value is
constant on the play. Along these strategies, the decision maker
ensures that when balancing between current payoff and future
states, he is not making irreversible mistakes, in the sense that
the uniform value does not decrease along the induced play.

In the third part we prove Theorem $\ref{theo1}(1)$. To prove the existence of a (non-pure) $0$-optimal strategy, we first show the existence of pure strategies such that the $\limsup$ of the long run expected average payoffs is the uniform value. Nevertheless the payoffs may not converge along the play induced by these strategies. We show that the decision maker
can choose a proper distribution over these strategies to ensure
that the expected average payoff converges.

The fourth part is dedicated to the proof of Theorem
$\ref{theo1}$(2). To construct a pure $0$-optimal strategy, instead of concatenating strategies one after the other, as is often done in the literature, we define a
sequence of strategies $(\sigma^l)_{l\geq 1}$ such that for every $l\geq 1$, $\sigma^l$ guarantees $v^*(x^l)-\varepsilon_l$ where $x^l\in X$ and $\varepsilon_l$ is a positive real number. We then split these strategies, seen as sequences of actions, into blocks of actions and construct a $0$-optimal strategy by playing these blocks in a proper order.
%
%In the following, we denote by $x_1$ the initial
%state.

\subsection{An example of a commutative deterministic MDP without $0$-optimal pure strategies}

In this section, we provide an example of a commutative
deterministic MDP with a uniform value in pure strategies that does
not have a pure $0$-optimal strategy. Before going into the details, we
outline the structure of the example. The set of states, which is
the countable set $\N\times \N$, is partitioned into countably many
sets, $\{h^0,h^1,...\}$, such that the payoff is constant on each
element of the partition; the payoff is $0$ on the set $h^0$ and
$1-\frac{1}{2^l}$ on  the set $h^l$, for all $l\geq 1$. We will
first check that for each $l\geq 1$, there exists a pure strategy
from the initial state $(0,0)$ that eventually belongs to the set
$h^l$. This will imply that the game starting at $(0,0)$ has a
uniform value equal to $1$. We will then prove that any $0$-optimal
pure strategy has to visit all sets $h^l$ and that when switching
from one set $h^i$ to another set $h^{i'}$, the induced play has to
spend many stages in the set $h^0$. The computation of the minimal
number of stages spent in the set $h^0$ shows that the average
expected payoff has to drop below $\frac{1}{2}$, which contradicts
the optimality of the strategy. Thus, there exists no $0$-optimal
pure strategy in the game starting at state $(0,0)$.

\begin{example}\label{escalier}
The set of states is $X=\N\times \N$ and there are only two actions
$R$ and $T$; the action $R$ increments the first coordinate and the action $T$ increments
the second one:
\begin{align*}
q((x,y),R)&=(x+1,y),\\
q((x,y),T)&=(x,y+1).
\end{align*}
Plainly the transition is deterministic and commutative.

For each $l\geq 1$, let $w^l=\sum_{m=1}^l \left(4^{m-1}-1 \right)=\frac{4^l-1}{3}-l$. We define the set $h^l \subset X$ by
\[h^l=\{(w_l,0) \} \cup \left\{(x,y), \ w^l+(y-1) \left(4^{l-1}-1 \right) \leq x \leq w^l+y \left(4^{l-1}-1 \right),\ x,y \geq 1 \right\}.\]
For example \[h^1=\{(0,y), \ y\geq 0\},\]
and
\[
h^2=\{(3,0) \} \cup \left(\cup_{y\geq 1}\{(3y,y),(3y+1,y),(3y+2,y),(3y+3,y)\} \right).
\]
For every $l\geq 1$, the set $h^l$ is the set of states obtained along the play induced by the sequence of actions $(T R^{4^{l-1}-1})^{\infty}$ from state $(w^l,0)$. We denote by $h^0$ the set of states not on any $h^l$, $l\geq 1$. Figure (\ref{contre-exemple}) shows the play associated to $h^1$, $h^2$, and $h^3$ with their respective payoffs. One can notice that the plays following these three sets separate from each other.

The payoff is $1-\frac{1}{2^l}$ in every state on the set $h^l$ and $0$ on the set $h^0:$
\[
g\left(x,y \right) =\begin{cases}1-\frac{1}{2^l} \text{ if }\  x \in \left[w^l+(y-1) \left(4^{l-1}-1 \right) ,w^l+y \left(4^{l-1}-1 \right)  \right] \\
0 \text{ otherwise.}
\end{cases}
\]

%This MDP is commutative and the transition is deterministic but there is no $0$-optimal pure strategy from $(0,0).$

\begin{figure}[h]
\begin{center}
\includegraphics[width=90mm, height=50mm]{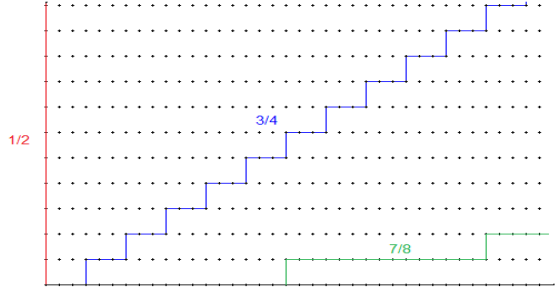}
  \caption{\label{contre-exemple} Payoff of the game on $h^1$,$h^2$ and $h^3$.}
\end{center}
\end{figure}
%
%The transition is clearly deterministic and it is immediate to check
%that it is commutative. Along a sequence of actions, the sequence of
%states visited has the shape of a stair. For $l=1$, the stair starts
%at $(0,0)$ and is degenerate since there is a unique infinite riser.
%For $l\geq 2$, the set of states on $h^l$ where the value is
%$1-\frac{1}{2^l}$ describes a stair. It begins at state $(0,w_l)$
%and it has a constant rise height of $1$ and a tread depth of
%$4^{l-1}-1$.
%Figure \ref{contre-exemple} shows the play for $l=1$, $l=2$ and $l=3$ with their respective payoffs. \\
\hspace{5mm}
%
%Formally, for all $z=(x,y)\in \N \times \N$ and $l\geq 1$, we say
%that $z$ is strictly in between $h^l$ and $h^{l+1}$, denoted by
%$z\in (h^l,h^{l+1})$, if and only if $w_l+(y-1)(4^{l-1}-1) <x <
%w_{l+1}+(y-1)(4^{l}-1)$. We denote by $[h^l,h^{l+1})$, the set where
%the left inequality is not strict. Let us call $h^0$ the set of
%states, which are strictly inbetween these stairs, such that
%$(h^l)_{l\geq 1}$ and $h^0$ induce a partition of $X$.\\

The uniform value exists in every state, is
equal to $1$, and the decision maker has $\varepsilon$-optimal pure
strategies: given an initial state, play $R$ until reaching some state in $h^l$ with $\frac{1}{2^l}\leq \varepsilon$ and then stay in $h^l$.\\

%We will then show that he has no $0$-optimal pure
%strategies. Given a play, we call the set of states in this history
%the path of the strategy and we say that two plays are crossing if
%there exists a common state on the two paths. For all $(p,q)\in \N
%\times \N$, there exists $l$ such that $(p,q)\in [h^l,h^{l+1})$.
%Therefore for all $l'\geq l+1$, there exists a finite number $n(l')$
%such that playing $n(l')$ times action $R$ leads to the path
%$h^{l'}$. Thus by following the path $h^{l'}$ after this stage, the
%decision maker can guarantee $1-\frac{1}{2^{l'}}$. The uniform value
%in pure strategies exists and is equal to $1$.
%% verstion courte

There exists no $0$-optimal pure strategy from
state $(0,0)$. Since there is only one player, the transition is
deterministic and the payoff depends only on the state, we can
identify a pure strategy with the sequence of states it selects on the play that it defines. Let $h=(z_1,...,z_t,...)$ be a $0$-optimal strategy. Since the play $h$ guarantees $1$, there
exists an increasing sequence of stages $(n_l)_{l\geq 1}$ such that $h$ crosses $h^l$ at stage $n_l$.

Let $m_l<n_{l+1}$ be the last time before $n_{l+1}$ such that $h$ intersects $h^l$. The reader can check that then $n_{l+1}-m_l \geq m_l$, and  every state of $h$ between stage $m_l+1$ and stage $n_{l+1}$ is in $h^0$. Therefore the expected average payoff at stage $n_{l+1}-1$ is less $\frac{1}{2} \times 1 + \frac{1}{2} \times 0=\frac{1}{2}$. \\

We now argue that there is a behavioral strategy that yields payoff $1$. We start by illustrating a strategy that yields an expected average payoff at least $\frac{3}{8}$ in all stages:
%
%\begin{example}
%\rm  We study the same game as in Example \ref{toto} and focus on a
%weaker property than the existence of a $0$-optimal strategy. The
%question is the existence of a strategy which starts at $(0,0)$,
%reaches the second stair $h^2$ and such that for every stage $n$,
%the mean-payoff is above $3/8$. It is impossible with pure
%strategies but possible if not playing purely.
%
%Any pure strategy, starting at $(0,0)$ and reaching the second
%stair, has to leave the first stair at some stage $n$. Then, by
%construction, the number of stages before reaching the second stair
%is more than $n$ and the mean expected payoff drops at $1/4$.
%
%
%Let us exhibit a strategy, such that the state reaches the second
%stair, but the payoff at each stage does not drop below $3/8$:
\begin{itemize}
\item with probability $1/4$, the decision maker plays $3$ times Right in
order to reach the set $h^2$ and then stays in $h^2$;

\item with probability $1/4$, the decision maker stays in the set $h^1$ for $3$ stages ($3$ times Top) then plays $9$ times Right in order to
reach the set $h^2$ and then stays in $h^2$;

\item with probability $1/4$, the decision maker stays in the set $h^1$
for $3+9=12$ stages then plays $36$ times Right in order to reach the set $h^2$ and then stays in $h^2$;

\item with probability $1/4$ the decision maker stays in the set $h^1$
for $3+9+36=48$ stages then plays $144$ times Right in order to reach the set $h^2$ and then stays in $h^2$.
\end{itemize}

Note that the state at stage $192$ is on $h^2$ and more precisely equal to $(48,144)$ whatever is the pure strategy chosen. The first strategy yields a payoff of $3/4$ except on the second and
third times the decision maker is playing right (stage
2-3), the second one yields a payoff of $1/2$ before stage $3$ (included) and a payoff of $3/4$ from stages $13$ (included), the third one yields a payoff of  $1/2$ before stage $12$ (included) and a payoff of $3/4$ from stage $49$ (included) and the fourth one yields a payoff of  $1/2$ before stage $48$ (included) and a payoff of $3/4$ from stages $193$ (included).

Thus for each stage $n$ up to stage $192$, there is at most one of the four pure strategies
that gives a daily payoff of $0$. The others strategies either
stay in $h^1$ or stay in $h^2$, and thus yield a stage
payoff at least $1/2$. Therefore, the expected payoff at each
stage is greater than $\frac{3}{4} * \frac{1}{2}=\frac{3}{8}$ and the expected average payoff until stage $n$ is greater
than $3/8$ for every $n\geq 1$. We managed to build a strategy going from $h^1$ to the set $h^2$ such that the expected average payoff does not drop below $\frac{3}{8}$.\\

We can iterate and switch from $h^2$ to $h^3$ without the payoff dropping below $\frac{7}{8}* \frac{3}{4}$ by considering $8$ different pure strategies. Repeating this procedure from $h^l$ to $h^{l+1}$ for every $l\geq 1$ will lead to a $0$-optimal strategy. To define properly a strategy which ensures an expected payoff $1$, we augment the strategy as in Section \ref{theo1.1}.

%By using more strategies, we could get a strategy such that the
%payoff does not drop below $1/2$-epsilon for all positive
%$\varepsilon$. In this example, in order to get a $0$-optimal strategy,
%we have to repeat the procedure between stair $2$ and stair $3$,
%stair $3$ and stair $4$ and to use at each switch more and more
%strategies to prevent the payoff for dropping too much.
\end{example}

\subsection{Existence of $\varepsilon$-optimal strategies with a constant value on the induced play}

In this part, we consider a commutative deterministic MDP with a
uniform value in pure strategies. We show that for all $x_1 \in X$
and all $\varepsilon>0$ there exists an $\varepsilon$-optimal pure
strategy in $\Gamma(x_1)$ such that the value is constant on the
induced play. Lehrer and Sorin \cite{Lehrer_92} showed that in
deterministic MDPs, given a sequence of actions, the value is always
non-increasing along the induced play. In particular, it is true
along the play induced by an $\varepsilon$-optimal pure strategy. We
need to define an $\varepsilon$-optimal pure strategy such that the
value
is non-decreasing.\\

To this end, we introduce a partial preorder on the
set of states such that, if $x'$ is greater than $x$, then $x'$ can be
\emph{reached} from $x$, i.e. there exists a finite sequence of actions such that $x'$ is on one play induced from $x$. Fix a state $x_1$ and let $x$ be a state which can be reached from $x_1$. By
commutativity, the order of actions is not relevant and we can represent the
state $x$ by a vector $m \in \N^I$, counting how many times each
action has to be played in order to reach $x$ from $x_1$. Let $M(x)$ be the set of all vectors representing the state $x$.

Given two vectors $m$ and $m'$ in $\R^I$, $m$ is greater than $m'$ if for every $i\in I$, $m(i) \geq m'(i)$. Given two states $x$ and $x'$, we say that $x$ is greater than $x'$, denoted $x\geq x'$, if there
exists $m \in M(x)$ and $m' \in M(x')$, such that $m \geq m'$. By construction, $x \geq x'$ implies that
$x$ can be reached from the state $x'$. Indeed if $x$ is
greater than $x'$, then there exists $m \in M(x)$ and $m' \in M(x')$
such that $m \geq m'$ in all coordinates. By playing
$(m-m')(i)$ times the action $i$ for every $i \in I$, the decision
maker can reach the state $x$ from $x'$.

\begin{lemma}\label{lemma43}
Consider a commutative deterministic MDP with a uniform value in
pure strategies. For all $x_1 \in X$ and all $\varepsilon>0$ there
exists an $\varepsilon$-optimal strategy in $\Gamma(x_1)$ such that
the value is non-decreasing, thus constant, on the induced play.
\end{lemma}

\noindent{\bf Proof:} Fix $x_1\in X$ and $\varepsilon>0$. We
construct a sequence of real numbers $(\varepsilon_l)_{l \geq 1}$ and a
sequence of strategies $(\sigma^l)_{l \geq 1}$ satisfying three
properties. For each $l\geq 1$, we denote by $(x^l_n)_{n\geq 1}$ the
sequence of states along $\sigma^l$. First, the sequence
$(\varepsilon_l)_{l \geq 1}$  is decreasing and $\varepsilon_1\leq
\varepsilon$ (property $(i)$). Second, for every $l\geq 1$ the strategy $\sigma^l$
is $\varepsilon_l$-optimal in the game $\Gamma(x_1)$ (property $(ii)$). Finally, given
any $l\geq 1$ and any stage $n\geq 1$, for every $l' \geq l$ there
exists a stage $n'$ such that $x_{n'}^{l'} \geq x_n^l$ (property $(iii)$). This
implies that $x_{n'}^{l'}$ is reachable from $x_n^l.$ Informally,
a decision maker who follows the strategy $\sigma^l$ can change his
mind in order to play better: at any stage he can stop following
$\sigma^l$, choose any $l'\geq l$, and play some actions such that
the play merges eventually with the play induced by $\sigma^{l'}$.
\newline

% plus precis

Let $(\varepsilon_l)_{l \geq 1}$ be a decreasing sequence of
positive numbers converging to $0$ such that
$\varepsilon_1=\varepsilon$. For each $l\geq 1$, let $\sigma^l$ be an $\varepsilon_l$-optimal pure strategy in
$\Gamma(x_1)$. We identify $\sigma^l$
with the sequence of actions $(i_1,i_2,...)$ it induces. By construction,
these sequences satisfy properties $(i)$ and $(ii)$. To satisfy property $(iii)$, we extract
a subsequence.

For all $l\geq 1$ and all $n \geq 1$, considering the strategy
$\sigma^l$ until stage $n$ defines a vector $m_n(\sigma^l)$ in $M(x^l_n)$. The sequence
$(m_n(\sigma^l))_{n \geq 1}$ is non-decreasing in every coordinate,
so we can define the limit vector $m_{\infty}(\sigma^l) \in (\N
\times \{\infty\})^{I}$.
By definition of the limit, for any $w\in \N^{I}$ such that
$w\leq m_{\infty}(\sigma^l)$, there
exists some stage $n$ such that $w \leq m_{n}(\sigma^l)$.

%This vector counts how many times
%each action is played. If one coordinate is equal to some integer
%$m$, the corresponding action is played $m$ times. If one coordinate
%is equal to $\infty$, the corresponding action is played infinitely
%often.

%
Since the number of actions is finite, we can choose a subsequence
of $(m_{\infty}(\sigma^l))_{l \geq 1}$ such that each coordinate is
non-decreasing in $l$. Informally, the closer to the value the decision maker
wants the payoff to be the more he has to play each
action. We keep the same notation, and denote by $(\varepsilon_l)_{l
\geq 1}$ and $(\sigma^l)_{l \geq 1}$ the sequences after extraction.

After extraction $\varepsilon_1$ is smaller than
$\varepsilon$. By definition, $\sigma^l$ is $\varepsilon_l$-optimal in the game
$\Gamma(x_1)$. Moreover, given two integers $l,l'$ such that $1\leq l \leq l'$, we have $m_\infty(\sigma^l) \leq
m_\infty(\sigma^{l'})$. Let $n$ be a positive integer, then
\[m_n(\sigma^l) \leq m_\infty(\sigma^l) \leq m_\infty(\sigma^{l'}).\]
By definition of $m_\infty(\sigma^{l'})$ as a limit, there exists a
stage $n'$ such that $m_{n'}(\sigma^{l'})$ is greater than
$m_n(\sigma^l)$, and thus $x^{l'}_{n'}$ is greater than $x^l_n$. The
subsequences $(\varepsilon_l)_{l \geq 1}$ and $(\sigma_l)_{l \geq 1}$ satisfy all the properties $(i)-(iii)$.\newline
%Let us check that the sequences $(\varepsilon_l)_{l \in \N}$ and
%$(\sigma^l)_{l \in \N}$ satisfy the required properties so the first property holds.

We now deduce that the value along $\sigma_1$ is non decreasing: for every $n\geq
1$, the uniform value in state $x^1_n$ is equal to the uniform value in the initial
state. Fix $n \geq 1$ and $l' \geq 1$. By construction, there exists $n'\geq n$ such that
$x_{n'}^{l'}$ can be reached from state $x_n^1$. Applying Lehrer and
Sorin \cite{Lehrer_92}, we know that the value is non increasing
along plays so $v^*(x^1_n)\geq v^*(x^{l'}_{n'}).$ Moreover, the strategy $\sigma^{l'}$ defines a continuation
strategy from $x^{l'}_{n'}$, which yields an average long-run payoff of at least $v^*(x_1)-\varepsilon_{l'}$. Thus, the uniform
value along the play induced by $\sigma_{l'}$ does not drop below
$v^*(x_1)-\varepsilon_{l'}$:
\[
v^*(x^{l'}_{n'})\geq v^*(x_1)-\varepsilon_{l'}.
\]
Considering both results together, we obtain that
\[
v^*(x^1_n)\geq v^*(x^{l'}_{n'}) \geq v^*(x_1)-\varepsilon_{l'}.
\]
Since it is true for every $l'\geq 1$, we deduce that the value is
non decreasing along $\sigma_1$. In order to conclude, notice that
$\varepsilon_1 \leq \varepsilon$, therefore $\sigma_1$ is
$\varepsilon$-optimal.\hfill $\Box$

\subsection{Proof of Theorem \ref{theo1}(1)}\label{theo1.1}

%\subsection{Existence of a $0$-optimal strategy in the general case (proof of Theorem \ref{theo1}(1))}

In this subsection, we prove  Theorem \ref{theo1}(1): in every commutative MDP with a
uniform value in pure strategies, there exists a $0$-optimal
strategy.

A strategy $\sigma$ is said to be \emph{partially $0$-optimal} if the limsup of the sequence of expected average payoffs is equal to the uniform value: $\limsup_n \gamma_n(x_1,\sigma)= v^*(x_1)$. We first deduce from Lemma \ref{lemma43} the existence of partially $0$-optimal pure strategies. As shown in Example
\ref{escalier}, expected average payoffs may not converge along partially $0$-optimal strategy and, in particular, can be small infinitely often. The key point of the proof of Theorem \ref{theo1}(1) is that
different partially $0$-optimal strategies have bad expected average payoff at different stages. By choosing a proper mixed strategy that is supported by pure partially $0$-optimal strategies, we can ensure that, at each stage, the probability to play one pure strategy with
a bad expected average  payoff is small.

%
%The average expected payoff along this strategy, get as close as wanted from the uniform value. Nevertheless,
%the payoff may also be bad infinitely often as shown in Example
%\ref{toto}. In this example, there are many partially $0$-optimal
%pure strategies but they all have infinitely often an average
%expected payoff smaller than $1/2$, far from the uniform value which
%is $1$.

We will first provide the formal definition of partially $0$-optimal strategies
and the concatenation of a sequence of strategies along a sequence
of stopping times. Then, we define two specific sequences such that the concatenated strategy  $\sigma^*$ is $0$-optimal. The proof of the
optimality of $\sigma^*$ is done in two steps: we check that the
support of $\sigma^*$ is included in the set of partially
$0$-optimal strategies, and that the probability to play a strategy
with a bad expected average payoff at stage $n$ converges to $0$ for $n$ sufficiently large.\\

We now start the proof of Theorem \ref{theo1}(1) by defining
formally a partially $0$-optimal strategy.
\begin{definition}
Let $\Gamma=(X,I,q,g)$ be an MDP and $v^*(x_1)$ be the uniform value of the MDP
starting at $x_1$. A strategy $\sigma$ is \emph{partially $0$-optimal} if
\[ \limsup_n \gamma_n(x_1,\sigma) = v^*(x_1).
\]
That is, for every $\varepsilon>0$, the long run expected average payoff is greater than $v^*(x_1)-\varepsilon$ infinitely often.
\end{definition}

We define the concatenation of strategies with respect to a sequence of stopping
times \footnote{A stopping time $u$ is a random variable such that the event $\{u \leq n \}$ is measurable with respect to the history up to stage $n$}. Let $(u_l)_{l \geq 2}$ be a sequence of increasing stopping times and
$(\sigma_l)_{l \geq 1}$ be a sequence of strategies. The concatenated strategy $\sigma^*$ is defined as follows. For every $t\geq 1$ and every $h_t=(x_1,i_1,j_1,...,x_t)$, let $l^*=l^*(h_t)=\sup\{ l, u_l(h_t) \leq t\}$ and $\sigma^*(h_t)=\sigma_{l^*}(h^{u_l^*}_t)$ where $h^{u_l^*}_t=(x_{u_l^*},i_{u_l^*},j_{u_l^*}...,x_t)$. Informally, for every $l\geq 2$, at stage $u_l$ the decision maker forgets the past history and follows $\sigma_{l}$.\\
% We start with the concatenation of two strategies. Given a
%stopping time $u$ and two strategies $\sigma$ and $\sigma'$, define the strategy
%$\sigma u \sigma'$ as follows: play $\sigma$ until $u$, then switch
%to $\sigma'$ as if the initial state is $x_u$ (and forget the history up to $u$). For every $t\geq 1$
%and every $h_t=(x_1,i_1,j_1,...,x_t)$, $(\sigma u
%\sigma')(h_t)=\sigma(h_t)$ if $u(h_t)>t$ and $(\sigma u
%\sigma')(h_t)=\sigma'(h^u_t)$ if $u(h_t)\leq t$ where
%$h^u_t=(x_u,i_u,j_u...,x_t)$. At stage $u$, the decision maker
%forgets the past and starts following $\sigma'$ as if it was the
%first stage of the game $\Gamma(x_u)$. The definition for a finite
%number of stopping time $u_2<...<u_m$ and a finite sequence of
%strategies $\sigma_1,...,\sigma_n$ is immediate.
%
%We now give a meaning to a countable concatenation. Let $(u_l)_{l
%\geq 2}$ be a sequence of increasing stopping times and
%$(\sigma_l)_{l \geq 1}$ be a sequence of strategies. For each $l\geq
%2$, we denote by $T_l$ the support of the stopping time $u_l$. For
%each $n$, we define $\sigma^*_n$ the concatenation associated to the
%finite number of stopping time $(u_2,...,u_n)$. The countable
%concatenation strategy $\sigma^*$ is the strategy which coincides
%with $\sigma^*_n$ before the random time $u_n$ for every $n\geq 1$.
%An increasing sequence of integers define naturally a sequence of
%increasing stopping time (deterministic) and thus a concatenated
%strategy. Note that the decision maker starts to play $\sigma_l$ at
%the random stage $u_l$.\newline

\underline{Definition of the $0$-optimal strategy:} Fix $x_1\in X$.
For every $t\geq 1$, we denote by $X(t)$ the set of states which can
be reached from $x_1$ in less than $t$ stages. Since the transition
is deterministic and the number of actions is finite, the set $X(t)$
is finite  for every $t\geq 1$. We choose two specific sequences of
stopping times and strategies and denote by $\sigma^*$ the
concatenation. Let $(\varepsilon_l)_{l \geq 1}$ be a decreasing
sequence of real numbers converging to $0$. For each $x\in X$ and
every integer $l\geq 1$, we denote by $\sigma_l(x)$ an
$\varepsilon_l$-optimal strategy in $\Gamma(x)$ such that the
uniform value is constant on the play, and let $N(l,x)$ be an
integer that satisfies
\begin{align}\label{tata}
  \forall n \geq N(l,x),\  \gamma_n(x,\sigma_l(x)) \geq v^*(x)-\varepsilon_l.
\end{align}
In any games longer than $N(l,x)$ stages, the average expected
payoff is close to the value, but the payoff in shorter games is not
controlled. The strategy $\sigma_l(x)$ exists by Lemma
\ref{lemma43}.

We now define the sequence of stopping times. For every $l\geq 1$, we define a set of stages $T_l$ and let $u_l$ be a stopping time uniformly distributed over $T_l$.
Start by setting $t_1=1$ and $T_1=\{1\}$. Let $l\geq
1$ and assume that the set $T_l$ is already defined. Denote $t_{l+1}=\left[
\frac{1}{\varepsilon_{l+1}}\right]+1$ and define the set
$T_{l+1}=\{T_{l+1}^{(1)},...,T_{l+1}^{(t_{l+1})}\}$ by induction:
\begin{align*}
T_{l+1}^{(1)} &=T_{l}^{(t_l)}+\max_{x\in X(T_{l}^{(t_l)})} N(l,x)+\left[\frac{1}{\varepsilon_l}+1 \right] T_{l}^{(t_l)},\\
T_{l+1}^{(2)}& =T_{l+1}^{(1)}+\max_{x \in X({T_{l+1}^{(1)}})} N(l+1,x),\\
...&, ...\\
T_{l+1}^{(t_{l+1})}&
=T_{l+1}^{(t_{l+1}-1)}+\max_{x \in X(T_{l+1}^{(t_{l+1}-1)})} N\left(l+1,x \right).
\end{align*}

Let $t\in T_{l+1}$, we call the smallest integer strictly greater than $t$ in $T_{l+1} \cup T_{l+2}$, the \emph{successor} of $t$. Formally, there exists $c_{l+1} \leq t_{l+1}$ such that $t=T_{l+1}^{(c_{l+1})}$. If $c_{l+1}$ is strictly smaller than $t_{l+1}$, the successor of $t$ is  $T_{l+1}^{(c_{l+1}+1)}$;
if $c_{l+1}=t_{l+1}$, then the successor of $t$ is $T_{l+2}^{(1)}$.\\

%Moreover the first possible realization of $u_l$ is big enough in
%order to outweigh everything that happened before.
%For each set $T_l$, we define the stopping time $u_l$ such that for
%all $m \in \{1,...,t_l\}$, $P(u_l=T_l^m)=\frac{1}{t_l}$ and
We make few comments on the definition of the set $T_{l+1}$. First, the number of stages between two different integers $t$ and $t'$ in $T_{l+1}$ is such that a strategy, which starts playing like
$\sigma_{l+1}(x_{t})$ at stage $t$ yields an expected
average payoff  between stage $t$ and stage $t'-1$ greater than $v^*(x_1)-\varepsilon_{l+1}$. Second, the weight of the first $T_{l}^{t_l}$ stages in a game of length $T_{l+1}^1$ is small.\\

%Given a stage $t$ in $T_l$, we call successor of $t$ either the next stage in $T_l$ if it exists or the first stage in $T_{l+1}$ otherwise.

%We denote by $\sigma^*$ the countable concatenation.

%\begin{remark}
%The strategy defined in the proof of Lemma \ref{partially} is the
%countable concatenation of $(\sigma_l)_{l\in \N}$ along the sequence
%of Dirac mass $(n_l)_{l\in \N}$.
%\end{remark}

We prove that the strategy $\sigma^*(x_1)$ is $0$-optimal. We consider here $\sigma^*(x_1)$ as a mixed strategy, i.e. a probability over pure strategies. More precisely, let $\Omega$ be the set of pure strategies defined as concatenation of a sequence of integers $(n_l)_{l\geq 2}$ with $n_l \in T_l$ for every $l\geq 2$ and the sequence of strategy $(\sigma_l)_{l\geq 1}$. $\sigma^*$ is a probability over $\Omega$.\\

We show that every pure strategy in $\Omega$ is partially $0$-optimal.

\begin{lemma}\label{partially}
Let $(n_l)_{l\geq 2}$ be a sequence of integers such that for every $l \geq 2$, $n_l \in T_l$. Denote by $\sigma$ the concatenated strategy induced by $(n_l)_{l \geq 2}$ and $(\sigma_l)_{l \geq 1}$.

The strategy $\sigma$ is partially $0$-optimal. Moreover, we have explicit lower bounds for specific stages. For every $l\geq 2$, let $n'_l$ be the successor of $n_l$. Then
\[
\forall l \geq 2, \ \forall n \in [n'_l-1,n_{l+1}-1], \
\gamma_n(x_1,\sigma) \geq v^*(x_1)- 2\varepsilon_{l-1}.
\]
\end{lemma}

%For every $l \geq 2$ let $n_l$ be a stage in
%$T_l$. There exists $c_l \leq t_l$ such that $n_l=T_l^{c_l}$. If
%$c_l$ is strictly smaller than $t_l$, set $n'_l=T_l^{c_l+1}$;
%if $c_l=t_l$, we define $n'_l=T_{l+1}^1$. By definition, $n_l < n'_l \leq n_{l+1}$.

%The payoff is not controlled inbeetween stage $n_l$ and $n'_l$, thus
%we consider the mean expected payoff stopping at stage $n'_l$.

\noindent{\bf Proof:} We first show that the sequence
$\gamma_{n_{l+1}-1}(x_1,\sigma)$ converges to the uniform value
$v^*(x_1)$ when $l$ goes to $\infty.$ At stage $n_l$, the strategy
$\sigma$ starts to follow an $\varepsilon_l$-optimal strategy from
the current state. By definition, $n'_l-n_l \geq
N(l,x_{n_l})$, and thus by Equation (\ref{tata})
\[
\gamma_{n_l,n'_l-1}(x_1,\sigma) = \gamma_{n'_l-n_l}(x_{n_l},\sigma_l(x_{n_l})) \geq
v^*(x_{n_l})-\varepsilon_l \geq v^*(x_1)-\varepsilon_l.
\]
More generally, for every $n \in [n'_l-1,n_{l+1}-1]$, we
have
\[
\gamma_{n_l,n}(x_1,\sigma) =
\gamma_{n-n_l+1}(x_{n_l},\sigma_l(x_{n_l})) \geq
v^*(x_1)-\varepsilon_l.
\]
In particular we have
\begin{align}\label{minor}
\gamma_{n_l,n_{l+1}-1}(x_1,\sigma) & \geq
v^*(x_1)-\varepsilon_l.
\end{align}

The expected average payoff between stage $n_l$ and
$n_{l+1}-1$ is greater than $v^*(x_1)-\varepsilon_l.$ It follows
that the sequence $(\gamma_{n_{l+1}-1}(x_1,\sigma))_{l \geq 1}$ converges to $v^*(x_1)$ and therefore the strategy $\sigma$ is partially $0$-optimal.\\

We now prove the second part of the lemma, giving some explicit
subsequences and bounds on the rate of convergence: for all $l \geq 2$, for
all $n$ between $n'_l-1$ and $n_{l+1}-1$, we prove that
\begin{align}\label{bouh}
\gamma_n(x_1,\sigma) \geq v^*(x_1)- 2\varepsilon_{l-1}.
\end{align}

Fix $l\geq 2$. We first prove this lower bound for the expected average payoff until stage $n_l-1$ (which is before $n'_l-1$).  By definition of $T_{l}^1$, the weight of the $n_{l-1}$ first stages is small in the MDP of length $n_{l}-1$:
\begin{align*}
\frac{n_{l-1}-1}{n_{l}-1} \leq \frac{n_{l-1}-1}{T_{l}^{(1)}-1} \leq &
\frac{T_{l-1}^{(t_{l-1})}}{T_{l-1}^{(t_{l-1})}+N(l-1,x_{T_{l-1}^{(t_{l-1})}})-1+\left[
\frac{1}{\varepsilon_{l-1}}+1 \right] T_{l-1}^{(t_{l-1})}} \\
& \leq \frac{T_{l-1}^{(t_{l-1})}}{\left[ \frac{1}{\varepsilon_{l-1}}+1
\right]T_{l-1}^{(t_{l-1})}} \leq \varepsilon_{l-1}.
\end{align*}

Using Equation (\ref{minor}) for $l'=l-1$ and the previous equation, it follows that
\begin{align*}
\gamma_{n_{l}-1}(x_1,\sigma) & = \left[ \frac{n_{l-1}-1}{n_{l}-1} \gamma_{n_{l-1}-1}(x_1,\sigma) + \frac{n_{l}-n_{l-1}}{n_{l}-1} \gamma_{n_{l-1}, n_{l}-1}(x_1,\sigma) \right], \\
 & \geq  \left[ \gamma_{n_{l-1},n_{l}-1}(x_1,\sigma) - \frac{n_{l-1}-1}{n_{l}-1} \right], \\
 & \geq v^*(x_1)-2\varepsilon_{l-1}.
\end{align*}

Let $n$ be a positive integer between $n'_l-1$ and $n_{l+1}-1$. The expected average payoff until stage $n$ is the convex combination of the expected average payoff until stage $n_l-1$ (before $n'_l-1$) and the average expected payoff between stages $n_l$ and $n$. Both of
these quantities are greater than $v^*(x_1)-2\varepsilon_{l-1}$, and therefore their convex
combination is greater than $v^*(x_1)-2\varepsilon_{l-1}$ as well. \hfill $\Box$

\begin{remark}
\rm Following the notation of Lemma \ref{partially}, if $n \in \cup_{ l\geq 2} [n_l,n'_l-2]$ then we only know that the $n$-stage expected average payoff is greater than $0$.
\end{remark}

\begin{lemma}
$\sigma^*(x_1)$ is a $0$-optimal strategy.
\end{lemma}

\noindent{\bf Proof:}
We consider here $\sigma^*(x_1)$ as a mixed strategy. Lemma \ref{partially} showed that with probability one the pure strategies in the support of $\sigma^*$ are partially $0$-optimal. \\

Let $l\geq 2$ and fix $n$ an integer in $[T^1_l,T^1_{l+1}-1]$.  We show that, with probability higher than
$1-\varepsilon_l$, the decision maker is following a pure strategy
giving an expected average payoff until stage $n$ higher than
$v^*(x_1)-2\varepsilon_{l-2}$.  \\

By definition, there exists a unique stage $n^*_l$ in $T_l$ such that $n$ is between  $n^*_l$ and $n'^*_l-1$ where $n'^*_l$ is the successor of $n^*_l$:
\begin{align}\label{firstblock}
n_l^* \leq n \leq n'^*_l-1.
\end{align}
  Let $\sigma$ be a pure strategy with positive probability under $\sigma^*$. There exists a sequence $(n_d)_{d\geq 2}$ such that for all $d\geq 2$, $n_d \in T_d$ and $\sigma$ is the concatenated strategy induced by $(n_d)_{d \geq 2}$ and $(\sigma_d)_{d \geq 1}$. We follow the previous notation and denote for every $d\geq 2$, the successor of $n_d$ by $n'_d$.  Since $n\in [T^1_l,T^1_{l+1}-1]$, by construction of the sets $T_{l-1}$,$T_l$, and $T_{l+1}$, we have
\begin{align}\label{secondblock}
n'_{l-1} \leq n \leq n_{l+1}-1.
\end{align}

\hspace{5mm}

We now use the inequalities (\ref{firstblock}) and
(\ref{secondblock}) to handle the three different cases depending on
the respective places of $n^*_l$, the beginning of the block
containing $n$, and $n_l$, the stage where the strategy $\sigma$ is
switching from an $\varepsilon_{l-1}$ strategy to an
$\varepsilon_l$-optimal strategy: $n_l> n^*_l $, $n_l< n^*_l$, and
$n_l=n^*_l$.\\

If $n_l> n^*_l $, then at stage $n$ the pure strategy $\sigma$ is still following the $\varepsilon_{l-1}$-optimal strategy from state $x_{n_{l-1}}$ and therefore yields a high expected average payoff.
Formally, we have $n'_{l-1}\leq n \leq n'^*_l-1 \leq n_l-1$, so that by Lemma \ref{partially} applied to $l'=l-1$,
\[
\gamma_n(x_1,\sigma) \geq v^*(x_1)-2\varepsilon_{l-2}.
\]

If $n_l< n^*_l $, then at stage $n$, the pure strategy $\sigma$ has already followed the $\varepsilon_{l}$-optimal strategy from state $x_{n_l}$ for a long time and thus yields a high expected average payoff. Formally, we have $n'_l \leq n^*_l \leq n \leq n_{l+1}-1$, so that by Lemma \ref{partially} applied to $l$,
\[
\gamma_n(x_1,\sigma) \geq v^*(x_1)-2\varepsilon_{l-1} \geq v^*(x_1)-2\varepsilon_{l-2}.
\]

Finally if $n_l=n^*_l$, we do not control the expected average payoff but by definition of the stopping time $u_l$ the probability of the event $\{n_l= n^*_l\}$ is smaller than $\varepsilon_l$ under $\sigma^*$.\\

We can now conclude. We denote by $\PP_{\sigma^*}$ the probability distribution induced by $\sigma^*$ on the set of pure strategy and $\E_{\sigma^*}$ the corresponding expectation. Since the payoffs are in $[0,1]$, it follows that
\[
\gamma_{n}(x_1,\sigma^*)= \E_{\sigma^*} \left( \gamma_{n}(x_1,\sigma) \right) \geq (1-\varepsilon_l) (v^*(x_1)-2\varepsilon_{l-2}) \geq v^*(x_1)-3\varepsilon_{l-2}.
\]
This is true for every $l \geq 1$ and every integer $n \in [T^1_l,T^1_{l+1}-1]$, therefore the expected average payoff converges to the uniform value: the strategy $\sigma^*$ is $0$-optimal.$\hfill \Box$

\subsection{Proof of Theorem \ref{theo1}(2)}

In this section, we prove Theorem \ref{theo1}(2): namely, if the set of states $X$ is a precompact metric space, the transition is $1$-Lipschitz, deterministic, and commutative, and the
payoff function is uniformly continuous, then there exists a pure
$0$-optimal strategy.

We will first justify the existence of the uniform value and that we can assume that the set of states is
compact. Then, we will define recursively a sequence of states
$(x^l)_{l\geq 1} $ such that $x^1=x_1$ and $x^{l+1}$ is a limit
point of states along an $\varepsilon_l$-optimal pure strategy
$\sigma^l(x^l)$ starting from $x^l$ where the value is constant on
the induced play. Therefore, the value in all these states is equal
to $v^*(x_1)$.
%All the evaluation functions have the same modulus of continuity, so all these states have the same value equal to $v^*(x_1)$.

For each $l\geq 1$ we will define by induction a sequence of stages $(n^l_k)_{k\geq 1}$ such that the sequence of states induced by $\sigma^l$
at stages $n^l_k$ converges to the limit point $x^{l+1}$. We impose
in addition conditions on $n^l_{l+1}$ and on the speed of convergence. This sequence
of stages splits the strategy $\sigma^l$ into a finite sequence of
streaks of actions. Given $k\geq 1$, we call an elementary block the
streak of actions played between stage $n^l_{k-1}$ and $n^l_{k}$.
Note that it has $n^l_{k}-n^l_{k-1}$ actions. By convention, the
first block starts at stage $n^l_0=1$.

We will define the $0$-optimal strategy $\sigma^*$ by playing these
elementary blocks in a specific order. The strategy $\sigma^*$ is
defined as a succession of two types of blocks $(A_l)_{l\geq 1}$ and
$(B_l)_{l\geq 1}$ such that for all $l\geq 1$, $A_l$ is composed of
$l+1$ consecutive elementary blocks from $\sigma^l(x^l)$ and
$B_l$ is composed of $l-1$ elementary blocks, one from each
$\sigma^{l'}(x^{l'})$ for $1\leq l' \leq l-1$:
\[\sigma^*=(A_1,B_1,A_2,B_2,A_3,....). \footnote{Recall that a pure strategy is identified with the sequence of actions it selects on the play path.}\]

Block $B_{l-1}$ ensures that the distance between the state at the beginning of block $A_l$ and $x^l$ is small. Block $A_l$ guarantees an expected average
payoff close to the value up to a function of
$\varepsilon_l$. Moreover, block $A_l$ is long enough for the total
expected average payoff of $\sigma^*$ at the end of $A_l$ to be
close to the value. It will follow that the strategy $\sigma^*$ is partially
$0$-optimal. The rest of the proof consists in showing that the
expected average payoff does not drop between these stages,
neither during block $B_{l+1}$ nor during the first stages of $A_{l+1}$.
It follows that the strategy is $0$-optimal.
\newline

%is always small compared to the total length before. We deduce that the mean average payoff in a finite game ending in $B_{l+1}$ or at the beginning of $A_{l+1}$ is close to the previous one. If we look at a game which does not end in the middle of $A_{l+1}$, the mean average payoff for the last stages is driven by $\sigma^l(x_l)$, so it is close to the value.

%\underline{Existence of a pure $0$-optimal strategy in the precompact case: }
%\noindent{\bf Proof:}[Theorem \ref{theo1} (ii)]

Let $\Gamma=(X,I,J,q,g)$ be a deterministic commutative MDP with a
precompact metric space, a uniformly continuous payoff function and
a $1$-Lipschitz transition. We first justify the existence of the
uniform value. We follow Section 6.1 in Renault \cite{Renault_2011}.
Let $\Psi=(Z,F,r)$ be an auxilliary dynamic programming problem. The
set of states is $Z=X\times I$, the correspondence is given by
\[
\forall (x,i)\in Z,\ F(x,i)=\{(q(x,a),a), \ a\in I\},
\]
and the payoff function is for all $(x,i)\in Z$, $r(x,i)=g(x,i).$ We consider on $Z$ the following metric $D((x,i),(x',i')=\max(d(x,x'),\delta_{i \neq i'}).$ The set of states $(Z,D)$ is precompact metric, the correspondence is $1$-Lipschtiz (i.e. for all $z,z' \in Z, z_1 \in F(z)$ there exists $z'_1 \in F(z)$ such that $D(z_1,z'_1) \leq D(z,z')$ ) and the payoff function is equicontinuous. By Corollary 3.9 of the same paper \cite{Renault_2011}, $\Psi$ has a uniform value for any initial state. We can deduce immediatly that $\Gamma(x_1)$ has a uniform value for every $x_1\in X$.\\

We now prove that we can assume that $X$ is compact. Define an MDP $\hat{\Gamma}(\hat{X},I,\hat{q},\hat{g})$ as follows: $\hat{X}$ is the Cauchy completion of $X$, $\hat{q}$ is the $1$-Lipschitz extension of $q$ to $\hat{X}$, $\hat{g}$ is the uniformly continuous extension\footnote{Note that an extension is not possible if the underlying function is only continuous.} of $g$ to $\hat{X}$. By Renault \cite{Renault_2011} and previous paragraph, both MDPs $\Gamma$ and $\hat{\Gamma}$ have a uniform value for any initial state.

Moreover the previous construction defines the transition on new states but does not change its value whenever it was already defined: for any state $x_1$ in $X$, $\hat{q}$ and $q$ coincides, as well as $g$ and $\hat{g}$. Therefore the MDPs $\Gamma(x_1)$ and $\hat{\Gamma}(x_1)$ are the same MDP on $X$. It follows that they have the same value. \\

%. This is done by considering an auxiliary MDP where the set of states is the completion of $X$ and by noting that the auxiliary MDP restricted to $X$ is identical to the original MDP. Informally, completing the state space in order to be compact does not
%change the MDP from $x_1$ if $x_1\in X$ since the transition and payoff are not changed on $X$.

%Note that if there exists a $0$-optimal pure strategy in the game
%$\hat{\Gamma}(x)$, then there exists a $0$-optimal pure strategy in
%the original game $\Gamma(x)$.

%
%Let $X$ be a precompact metric space, $I$ be a finite
%set of actions, $q$ be a $1$-Lipschitz transition and $g$ be a
%uniformly continuous payoff function. We denote by
%$\Gamma=(X,I,q,g)$ this MDP.
% We consider the game
%$\hat{\Gamma}=(\hat{X},I,\hat{q},\hat{g})$.
%
%For any $x_1$ in the original set $X$, the strategies in both games generate the same set
%of probabilities over histories. Indeed, if $x_1$ is an initial
%state in $X$, any play in the game $\hat{\Gamma}$ from $x_1$ stays
%in $X$ by an immediate induction: $\hat{q}$ and $q$ coincides on $X$
%and $q$ has values in $\Delta_f(X)$, so starting in $X$, the state
%at the next stage is with probability $1$ in $X$. Therefore, given a
%state $x_1\in X$, the uniform value are equal at $x_1$ and a
%$0$-optimal strategy in $\hat{\Gamma}(x_1)$ is well defined and
%$0$-optimal in the game $\Gamma(x_1)$. \newline

In the following we assume that $X$ is compact. Let $x_1 \in X$ and let
$(\varepsilon_l)_{l \geq 1}$ be a decreasing sequence of positive real
numbers that converges to $0$. For each $x\in X$ and $l\geq 1$ denote by $\sigma_l(x)$ an $\varepsilon_l$-optimal pure strategy in
$\Gamma(x)$ such that the value along the induced play is constant, and by $N(l,x)$ an
integer such that
\[ \forall n \geq N(l,x),\  \gamma_n(x,\sigma_l(x)) \geq v^*(x)-\varepsilon_l.\]
%
%Given a sequence of actions $(i_n)_{n\in \N}$, we denote by
%$(x_t)_{t\geq 1}$ the sequence obtained from $x$ by playing $\sigma$
%and $(x'_t)_{t\in \geq 1}$ the sequence obtained from $x'$ by
%playing $\sigma$.

\noindent Since $g$ is uniformly continuous, there exists
$(\eta_l)_{l\geq 1}$ such that
\[\forall x,x' \in X,\ d(x,x')\leq \eta_l, \ \forall a\in \Delta(I), \ | g(x,a)-g(x',a)| \leq \varepsilon_l.\]

\noindent Let $\sigma=(i_t)_{t\geq 1} \in I^\infty$ be an infinite sequence of actions and let $x_1$ and $x_1'$ be two initial states. For every $n\geq 1$, the distance between $x_n$, the state at stage $n$ obtained along the play induced by $x_1$ and $\sigma$, and $x'_n$, the state at stage $n$ obtained along the play induced by $x'_1$ and $\sigma$, is smaller than $d(x_1,x'_1)$. It follows that
\[
\forall x_1,x'_1 \in X,\ s.t. \  d(x_1,x'_1)\leq \eta_l, \  \forall \sigma=(i_t)_{t\geq 1} \in I^\infty,\ \forall n\geq 1, \ |\gamma_n(x_1,\sigma)-\gamma_n(x'_1,\sigma)| \leq \varepsilon_l.
\]

\underline{Definition of the strategy $\sigma^*$:} Let $x^1=x_1$. Given
$(x^j)_{1\leq j\leq l}$ define $x^{l+1}$ to be a limit point of the
play $(x^l,\sigma_l(x^l))$. Since the value is constant on the play
induced by $(x^l,\sigma_l(x^l))$, the uniform value in $x^{l+1}$ is also equal
to $v^*(x_1)$. To construct the $0$-optimal strategy, we split each
play $\sigma_j(x^j)$ into blocks by induction on $j$.

Let us assume
that $(n^j_{k})_{k\geq 1}$ have been defined for every $j\leq l-1$, i.e. the splittings of all strategies $\{\sigma_1(x^1),...,\sigma_{l-1}(x^{l-1}\}$ have been defined. We now split the sequence $\sigma_l(x^l)$.

Define $L_l=1+\sum_{j\leq {l-1}} (n^j_l-1)$, which depends only on
the sequences for $j\leq (l-1)$. We denote by $(x^l_{n})_{n\geq 1}$
the sequence of states along $(x^l,\sigma_l(x^l))$. Let us define
the sequence of stages $(n^l_{k})_{k\geq 1}$ such that it satisfies four
properties. The three first properties are restriction on $n^l_{l+1}$ and the last one is a restriction on the rate of convergence to $x^{l+1}$. First the strategy $\sigma_l(x^l)$
guarantees in $\Gamma(x^l)$ the value with an error less than
$\varepsilon_l$ in all games longer than $n^l_{l+1}$:
\begin{eqnarray}\label{eq:1}
n^l_{l+1} &\geq N(l,x^l).
\end{eqnarray}
Second, $L_l$ is small compared to $n^l_{l+1}$:
\begin{eqnarray}\label{eq:2}
\frac{L_l}{n^l_{l+1}} \leq \varepsilon_l.
\end{eqnarray}
Third,
\begin{eqnarray} \label{eq:4}
\frac{N(l+1,x^{l+1})+\sum_{j=1}^{l-1}{\left( n^{j}_{l+1}-n^j_l
\right)}}{n^l_{l+1}} \leq \varepsilon_l.
\end{eqnarray}
Finally, at the beginning of the $k$-th block of this decomposition
the state is close to the limit point
\begin{eqnarray} \label{eq:3}
d(x^l_{n^l_k},x^{l+1}) \leq \frac{\eta_{k}}{k-1}.
\end{eqnarray}

%Note that $L_l$ is the length given by the concatenation the $n^{l'}_l$ first stage for each $l'\leq l$.
%Thus $n^l_{l+1}$ is big compared to the total length of the $l-1$ first intervals of the block $l$ and the time such that the mean-payoff is optimal if we start to play $\sigma_{l+1}$ in $\Gamma(x^{l+1})$.
\noindent Fix $l\geq 1$. We define $A_l$ to be the finite sequence of actions given by $\sigma^l(x^l)$ between stage $1$ and stage $n^l_{l+1}$. In term of elementary blocks, it is composed of the first $l+1$ elementary blocks of $\sigma^l(x^l)$ and is composed of $n^l_{l+1}-1$ actions. We define $B_l$ as the sequence of actions where the decision maker is playing, for each $l'<l$, the elementary block of $\sigma^{l'}(x^{l'})$ between stages $n^{l'}_{l}$ and $n^{l'}_{l+1}$. Thus $B_l$ is the concatenation of $l-1$ elementary blocks. Moreover the number of actions in $B_l$ is $b_l=\sum_{j=1}^{l-1}{\left( n^{j}_{l+1}-n^j_l \right)}$, which appeared in (\ref{eq:4}). The strategy $\sigma^*$ is the sequence of actions given by the alternating sequence $(A_l,B_l)_{l\geq 1}$.\\

We now show that the strategy $ \sigma^*$ is $0$-optimal.\\

%We first prove that $\sigma^*$ is partially $0$-optimal. For each $l\geq 1$, we show that the state
%at the beginning of $A_l$ is close to $x^l$ and we deduce that the
%expected average payoff of $\sigma^*$ at the end of $A_l$ is bigger
%than $v^*(x_1)- 3\varepsilon_l$. Thus $\sigma^*$ is partially
%$0$-optimal.
We first prove that the state at the beginning of $A_l$ is close to
$x^l$. Therefore the expected average payoff of $\sigma^*$ at the
end of $A_l$ is bigger than $v^*(x_1)- 3\varepsilon_l$ and
$\sigma^*$ is partially $0$-optimal.

\begin{lemma}\label{partopti}
The payoff at the end of $A_{l}$ is greater than
$v^*(x_1)-3\varepsilon_l$:
\[\gamma_{L_l+ n^l_{l+1}-1}(x_1,\sigma^*) \geq v^*(x_1)-3\varepsilon_l\]
\end{lemma}

\begin{corollary}
The strategy $\sigma^*$ is partially $0$-optimal.
\end{corollary}

\noindent{\bf Proof of Lemma \ref{partopti}:} Let us denote by $(x_n)_{n\geq 1}$ the
sequence of states on the play induced by $\sigma^*$.

We first prove that the state at the beginning of $A_l$ is close to $x^l$. One can verify that the
first stage of $A_l$ is the stage $L_l=1+\sum_{j\leq {l-1}}
(n^j_l-1)$. By definition, at stage $L_l$ for each $l'\leq l-1$, all
first $l$ elementary blocks of $\sigma^{l'}(x^{l'})$ have been
played: all of the first $l'+1$ on block $A_{l'}$ and then one after each other
 in the blocks $B_j$ for $j\in [l'+1,l-1]$. By commutativity, the state does not
depend on the order of actions and the state is the same as after the sequence
$\sigma'$ where the decision maker plays $\sigma_1(x^1)$ for
$n^1_{l}-1$ stages, $\sigma_2(x^2)$ for $n^2_{l}-1$ stages,..., and 
$\sigma_{l-1}(x^{l-1})$ for $n^{l-1}_{l}-1$ stages.

For each strategy $\sigma_j$, Equation (\ref{eq:3}) implies that the distance
between $x^{j+1}$ and the state at stage $n^j_l$ on the play from $x^j$ is less than $\frac{\eta_{l}}{l-1}$
for each $j\in \{1,...,l-1\}$. The map $q$ is $1$-Lipschitz, so the
distances sum up and an immediate induction implies that
\begin{align}\label{eq:5}
d(x_{L_l},x^l)\leq \eta_{l}.
\end{align}

Let us now compute the payoff in the MDP of length $L_l+n^l_{l+1}-1$, i.e. until the end of $A_l$. Equation (\ref{eq:2}) ensures that the
payoff is almost equal to the payoff between stages $L_l$ and
$L_l+n^l_{l+1}-1$:
\begin{align*}
\gamma_{L_l+n^l_{l+1}-1}(x_1,\sigma^*) & = \frac{L_l-1}{L_l+n^l_{l+1}-1}\gamma_{L_l-1}(x_1, \sigma^*) + \frac{n^l_{l+1}}{L_l+n^l_{l+1}-1}\gamma_{L_l,L_l+n^l_{l+1}-1}(x_1, \sigma^*) \\
& \geq \frac{n^l_{l+1}}{L_l+n^l_{l+1}-1}\gamma_{L_l,L_l+n^l_{l+1}-1}(x_1, \sigma^*) \\
& \geq \gamma_{L_l,L_l+n^l_{l+1}-1}(x_1, \sigma^*)- \frac{L_l-1}{L_l+n^l_{l+1}-1} \\
& \geq \gamma_{L_l,L_l+n^l_{l+1}-1}(x_1, \sigma^*)- \varepsilon_l.
\end{align*}
Moreover $\sigma^*$ plays like an $\varepsilon_l$-optimal strategy in $\Gamma(x^{l})$ between stages $L_l$ and $L_l+n^l_{l+1}-1$, and the distance between $x_{L_l}$ and $x^l$ is less than $\eta_l$ by Equation (\ref{eq:5}). Therefore, by Equation (\ref{eq:1}) we have
\begin{align*}
\gamma_{L_l+n^l_{l+1}-1}(x_1,\sigma^*)  & \geq \gamma_{n^l_{l+1}}(x_{L_{l}}, \sigma_l(x^{l}))- \varepsilon_l\\
& \geq \gamma_{n^l_{l+1}}(x^l, \sigma_l(x^{l}))- 2 \varepsilon_l\\
& \geq v^*(x_1)- 3\varepsilon_l.
\end{align*}$\hfill \Box$ \\

%Let show that in fact the payoff is controlled on the rest of the play.
%For each $l\geq 1$, the expected average payoff is good at the end
%of the block $A_l$.

We now check that the average expected payoff does not
drop between these stages. We distinguish between two different cases: if $n\in [L_l+n^l_{l+1}-1,
L_{l+1}+N(l+1,x^{l+1})]$ or if $n\in [L_{l+1}+N(l+1,x^{l+1}),L_{l+1}+n^{l+1}_{l+2}-1]$.

In the first case, the MDP ends at a stage in $B_{l}$ or in the beginning of block $A_{l+1}$. Equation (\ref{eq:4}) implies that the length of the game is almost equal to $L_l+n^l_{l+1}-1$, therefore the expected average payoff is greater than $v^*(x_1)- 4\varepsilon_l$.

In the second case, the MDP ends in the middle of block $A_{l+1}$. The expected average payoff is the convex combination of the expected average payoff until $L_{l+1}-1$ and the average expected payoff between $L_{l+1}$ and $n$. We check that both of them are high and we deduce that the expected average payoff is greater than $v^*(x_1)- 4\varepsilon_l$.

%
%We now prove that the payoff does not drop until the
%strategy played in $A_{l+1}$ ensures a good payoff itself.

\begin{lemma}\label{control1}
Let $n\in [L_l+n^l_{l+1}-1,
L_{l+1}+N(l+1,x^{l+1})]$. Then
\[
\gamma_n(x_1,\sigma^*) \geq v^*(x_1)-4\varepsilon_l.
\]
The expected average payoff in any $n$-stage MDP such that $n$ is in the middle of block $B_l$ or at the beginning of block $A_{l+1}$ is greater than $v^*(x_1)-4\varepsilon_l.$
\end{lemma}

%First we focus on the case where we consider a stage in $B_l$ or at less then $N(l+1,x^{l+1})$ after the beginning of $A_{l+1}$.

\noindent{\bf Proof:} The key point is that the number of stages is
close to the case of Lemma \ref{partopti}. Let $n\in [L_l+n^l_{l+1}-1,
L_{l+1}+N(l+1,x^{l+1})]$. By equation (\ref{eq:4}), we have
\begin{align*}
n-L_l-n^l_{l+1}+1 &\leq N(l+1,x^{l+1})+\sum_{j=1}^{l-1}{(n^{j}_{l+1}-n^j_l)} \\
                &\leq \varepsilon_l n_{l+1}^l.
\end{align*}

\noindent It follows that
\begin{align*}
\gamma_{n}(x_1,\sigma^*) & = \frac{L_l+n^l_{l+1}-1}{n}\gamma_{L_l+n^l_{l+1}-1}(x_1, \sigma^*) + \frac{n-L_l-n^l_{l+1}+1}{n}\gamma_{L_l+n^l_{l+1},n}(x_1, \sigma^*) \\
& \geq \frac{L_l+n^l_{l+1}-1}{n}\gamma_{L_l+n^l_{l+1}-1}(x_1, \sigma^*) \\
& \geq \gamma_{L_l+n^l_{l+1}-1}(x_1, \sigma^*)- \frac{n-L_l-n^l_{l+1}+1}{n} \\
& \geq v^*(x_1)- 3\varepsilon_l - \frac{n-L_l-n^l_{l+1}+1}{n^l_{l+1}} \\
& \geq v^*(x_1)- 4\varepsilon_l. \hfill \Box
\end{align*}

\begin{lemma}\label{control2}
Let $n\in [L_{l+1}+N(l+1,x^{l+1}),L_{l+1}+n^{l+1}_{l+2}-1]$. Then
\[
\gamma_n(x_1,\sigma^*) \geq v^*(x_1)-4\varepsilon_l.
\]
The payoff in any $n$-stage MDP stopping in the middle of block $A_{l+1}$ is greater than $v^*(x_1)-4\varepsilon_l.$
\end{lemma}

\noindent{\bf Proof:} Let $n\in [L_{l+1}+N(l+1,x^{l+1}),L_{l+1}+n^{l+1}_{l+2}-1]$. The expected  average payoff is the convex combination  of the expected average payoff until $L_l+n^l_{l+1}-1$ and the expected average payoff between $L_l+n^l_{l+1}-1$ and $n$. It follows that
\begin{align*}
\gamma_{n}(x_1,\sigma^*) & = \frac{L_{l+1}-1}{n}\gamma_{L_{l+1}-1}(x_1, \sigma^*) + \frac{n-(L_{l+1}-1)}{n}\gamma_{L_{l+1},n}(x_1, \sigma^*) \\
& = \frac{L_{l+1}-1}{n}\gamma_{L_{l+1}-1}(x_1, \sigma^*) + \frac{n-(L_{l+1}-1)}{n}\gamma_{n-L_{l+1}+1}(x_{L_{l+1}}, \sigma_{l+1}(x_{L_{l+1}}))\\
& \geq \frac{L_{l+1}-1}{n}(v^*(x_1)-4 \varepsilon_l) + \frac{n-(L_{l+1}-1)}{n}(v^*(x^{l+1})-2\varepsilon_{l+1}) \\
& \geq v^*(x_1)- 4\varepsilon_l.
\end{align*}

\noindent The expected average payoff is greater than $v^*(x_1)-4\varepsilon_l$. \hfill $\Box$\\

Lemma \ref{control1} and Lemma \ref{control2} are true for every $l\geq 1$, therefore the strategy $\sigma^*$ is pure and $0$-optimal at $x_1$, which concludes the proof.

\section{Commutative stochastic games.}\label{comm}

In this section, we focus on commutative stochastic games and state-blind repeated games. In Section \ref{absorb}, we show that the class of absorbing games is in fact a subclass of commutative stochastic games.  We show that each absorbing state can be replaced by a non-absorbing state leading to some new states, which are useless from a
strategic point of view but designed in order to fulfill the
commutativity assumption. In Section \ref{preuvecomm}, we prove the existence of the
uniform value in stochastic games with a deterministic commutative
$1$-Lipschitz transition (Theorem \ref{theo2}). In Section \ref{preuveblindstate}, we deduce the existence of the uniform value in state blind commutative repeated games (Corollary \ref{mdp_noir}). In Section \ref{extension}, we provide some generalizations.

%We show that we can introduce finite stochastic games and solve the problem by induction on the initial state and by using the result of Mertens and Neyman \cite{Mertens_81} for finite stochastic games.

\subsection{Absorbing games}\label{absorb}

Absorbing games were introduced by Kohlberg \cite{Kohlberg_74a}.
They are stochastic games with a single non-absorbing state. An
absorbing game is thus given by $\Gamma=(\{\alpha\} \cup X,I,J,q,g)$
where $\alpha$ is the unique non-absorbing state and all states
$x\in X$ are absorbing: $q(x,i,j)(x)=1\ \forall x\in X, i\in I, j\in
J.$ The state $\alpha$ is the only state where the players have an
influence on the payoff and on future states. For each action pair
$(i,j)\in I\times J$, we denote by $q(\alpha,i,j)(X)$ the total
probability to reach an absorbing state by playing the action pair
$(i,j)$.

\begin{proposition} \label{absorbant}
Let $\Gamma=(\{\alpha\} \cup X,I,J,q,g)$ be an absorbing game. There exists a commutative game $\Gamma'=(X',I,J,q',g')$ and a state
$\alpha_2' \in X'$ such that for all $n\geq 1$,
$v_n(\alpha)=v'_n(\alpha_2').$ Moreover a player can guarantee
$w$ in $\Gamma'(\alpha_2')$ if and only if he can guarantee
$w$ in $\Gamma(\alpha)$.
\end{proposition}

\noindent{\bf Proof:}
Let $q(\alpha,i,j|X)$ be the conditional probability on $X$ if the action pair $(i,j)$ is played and there has
been absorption. Define an auxiliary commmutative game $\Gamma'=(X',I',J',q',g')$ as follows. The action spaces are $I'=I$ and $J'=J$. For each $i\in I$ (resp. $j\in J$), we define a
new state $x_i$ (resp. $x_j$). The state space is given by
$X'=X_{I} \times X_{J}$, where $X_I=\{\alpha'\} \cup\{x_i ,\ i\in I \} \cup\{ \omega\}$ and $X_J=\{\alpha' \} \cup\{x_j, \ j\in J \} \cup\{ \omega\}$. In the following, we denote $(\alpha',\alpha')$ by $\alpha'_2$. The payoff function is defined by
\begin{center}
\begin{tabular}{lr@{}ll}
$\forall i,i' \in I, \forall j,j' \in J,$ & $g'(\alpha'_2,i,j)$ & $=g(\alpha,i,j)$, \\
                   & $g'((x_{i'},x_{j'}),i,j)$ & $=\E_{q(\alpha,i',j'|X)}(g(x))$,\\
                   & $g'((x_{i'},\omega),i,j)$ & $=1$,&\\
                   & $g'((\omega,x_{j'}),i,j)$ & $=0$,&\\
                   & $g'((\omega,\omega),i,j)$ & $=1/2$.&
\end{tabular}
\end{center}

The payoff function in $\Gamma'$ reflects the role of the
different states. The state
$\alpha'_2$ is a substitute of the state $\alpha$, and for
each pair $(i',j')$, the state $(x_{i'},x_{j'})$ replaces the
absorption occurring in state $\alpha$ by playing the action pair
$(i',j')$. This state will not be absorbing but an equilibrium at
$(x_{i'},x_{j'})$ is to stay in this state. If player $1$ deviates, then
with some probability the state will remain $(x_{i'},x_{j'})$ and
with the remaining probability the new state will be $(\omega,x_{j'})$, where
player $2$ can guarantee a payoff of $0$. Similarly, if player $2$ deviates,
then the new state will remain $(x_{i'},x_{j'})$ with some probability
and with the remaining probability it will be $(x_{i'},\omega)$ where player
$1$ can guarantee a payoff of $1$.\\

The transition $q'$ is defined in three steps: we define two
transitions $s_I$ on $X_{I}$ controlled only by player $1$
and $s_J$ on $X_{J}$ controlled only by player $2$. We then
consider the product transition corresponding to the absorbing
part of $q$, and finally we define $q'$. At each step, we check that
the transition is commutative. We define $s_{I}$ and $s_{J}$ by
\begin{center}
\begin{tabular}{rr@{}lrr@{}l}
$\forall i,i' \in I,$ & $s_{I}(\alpha',i)=$ & $\ x_{i}$,  & $\ \forall j,j' \in J,$& $s_{J}(\alpha',j)=$ & $\ x_{j}$,\\

                    & $s_{I}(x_{i'},i)=$ & $\begin{cases}
                      x_{i'} & \text{ if } i=i',\\
                      \omega & \text{ if } i \neq i',\\
                     \end{cases}$      &                                          & $s_{J}(x_{j'},j)=$ & $\begin{cases}
                                                                                               x_{j'} & \text{ if } j=j',\\
                                                                                               \omega & \text{ if } j \neq j', \\
                                                                                                 \end{cases}$\\

                   & $s_{I}(\omega,i)=$ & $\omega$,  &                 & $s_{J}(\omega,j)=$ & $\omega$.

\end{tabular}
\end{center}
%\begin{center}
%\begin{tabular}{}
%\end{tabular}
%\end{center}

\noindent We now verify that $s_I$ is commutative. A similar argument shows that $s_J$ is commutative. Let $i$ and $i'$ be two actions of player $1$. It is sufficient to check that $s_I$ commutes when $i \neq i'$. However, if player $1$ plays $i$ and $i'$, the state after two stages is $\omega$ regardless  of the initial state and of the order in which he plays these actions.

\begin{figure}
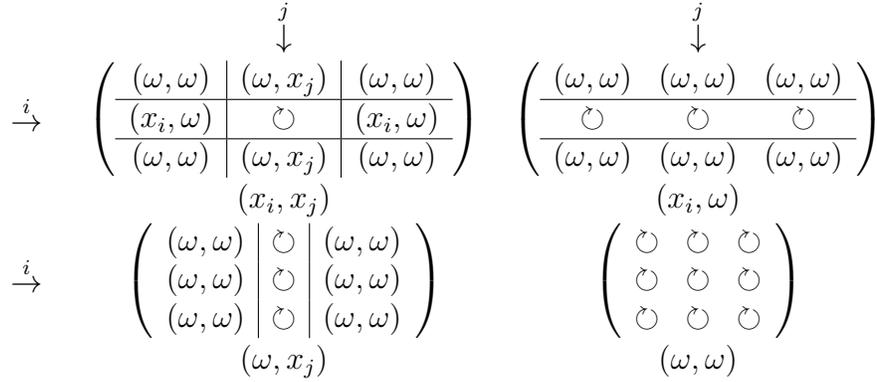

\begin{center}
\begin{tabular}{ccc}
 & $\begin{array}{ccc}
\ & \overset{j}{\downarrow} & \ \\
\end{array}$& $\begin{array}{ccc}
\ & \overset{j}{\downarrow} & \ \\
\end{array}$ \\
$\begin{array}{c}
\ \\
\overset{i}{\rightarrow}\\
\ \\
\end{array}$ & $\left( \begin{array}{c|c|c}
(\omega,\omega) & (\omega,x_j) & (\omega,\omega) \\
\hline
(x_i,\omega) & \circlearrowright & (x_i,\omega) \\
\hline
(\omega,\omega) & (\omega,x_j) & (\omega,\omega) \\
\end{array} \right)$ &

$ \left( \begin{array}{ccc}
(\omega,\omega) & (\omega,\omega) & (\omega,\omega) \\
\hline
\circlearrowright & \circlearrowright & \circlearrowright \\
\hline
(\omega,\omega) & (\omega,\omega) & (\omega,\omega) \\
\end{array} \right)$\\
& $(x_i,x_j)$ & $(x_i,\omega)$ \\

$\begin{array}{c}
\ \\
\overset{i}{\rightarrow}\\
\ \\
\end{array}$ & $\left( \begin{array}{c|c|c}
(\omega,\omega) & \circlearrowright & (\omega,\omega) \\
(\omega,\omega) & \circlearrowright & (\omega,\omega) \\
(\omega,\omega) & \circlearrowright & (\omega,\omega) \\
\end{array} \right)$ & $\left( \begin{array}{ccc}
\circlearrowright & \circlearrowright & \circlearrowright \\
\circlearrowright & \circlearrowright & \circlearrowright \\
\circlearrowright & \circlearrowright & \circlearrowright \\
\end{array} \right)$ \\
& $(\omega,x_j)$& $(\omega,\omega)$
\end{tabular}
\end{center}
 \caption{\label{abs_fig} A graphic depiction of $s$.}
\end{figure}

Let $s$ be the transition on  $X_{I} \times X_{J}$ defined by $s((x,y),(i,j))=(s_I(x,i),s_J(y,j))$. The reader can verify that $s$ is commutative; it is depicted graphically in Figure \ref{abs_fig}.

Let $q'$ be defined as follows: $q'(x,i,j)=q(\alpha,i,j)(\alpha)\delta_{x} + q(\alpha,i,j)(X)\delta_{s(x,i,j)}$ for all $x \in X'$, for all $i \in I$, and for all $j \in J$.
Thus for all $x\in X$, for all $i,i' \in I$ and for all $j,j' \in J$ we have
\begin{align}\label{eq_comm}
%q(q(x,i,j),i',j') &=q\left((1-g(i,j))\delta_{x} + g(i,j)\delta_{s(x,i,j)},i',j' \right) \\
%                  &=(1-g(i,j)q(x,i',j') + g(i,j) q(s(x,i,j),i',j') \\
\begin{split}
\widetilde{q}'(q'(x,i,j)&,i',j')=q(\alpha,i,j)(\alpha)q(\alpha,i',j')(\alpha)\delta_x+q(\alpha,i,j)(\alpha)q(\alpha,i',j')(X)\delta_{s(x,i',j')} \\
                   &+q(\alpha,i,j)(X)q(\alpha,i',j')(\alpha)\delta_{s(x,i,j)}+q(\alpha,i,j)(X)q(\alpha,i',j')(X) \delta_{s(s(x,i,j),i',j')}.
\end{split}
\end{align}
The right hand side of Equation (\ref{eq_comm}) is symmetric between $(i,j)$ and $(i',j')$ except the last term that involves $s$. Since $s$ is commutative, so is $q'$. Note that $\widetilde{q}$ may not be the product of one function depending on $I$ and one function depending on $J$. \newline

Fix $n\geq 1$. We prove that the $n$-stage values of $\Gamma(\alpha)$ and the $n$-stage values of $\Gamma'(\alpha'_2)$ are equal. Since the state $(\omega,\omega)$ is absorbing, the value is equal to $1/2$, the stage payoff. For all $i'$ in
$I$, the state $(x_{i'},\omega)$ is controlled by player $1$. His optimal action is $i'$ which guarantees him a payoff of $1$. The situation is symmetric for $(\omega,x_{j'})$, so for all $j'\in J$,
$v'_n((\omega,x_{j'}))=0$. Fix $(i',j') \in I\times J$. The action $i'$(resp. $j'$) is
optimal for player $1$ (resp. $2$) in state $(x_{i'},x_{j'})$ thus $v'_n(x_{i',j'})=\E_{q(\alpha,i',j'|X)}(g(x))$. The stage payoffs and the continuation values are equal in both the game $\Gamma(\alpha)$ and the game $\Gamma'(\alpha_2)$ so the values
in $\alpha$ and in $\alpha'_2$ are equal.\\

Finally there is a correspondence between strategies. Given a strategy $\sigma$ for player $1$ in the absorbing game $\Gamma$ that guarantees $w$, define $\sigma'$ in $\Gamma'$ by $\sigma'(\alpha'_2)=\sigma(\alpha)$ and for all $i'\in I$, $\sigma'(x_{i'},.)=i'$. For all $i'\in I$ and $j'\in J$, this strategy guarantees the payoff $\E_{q(\alpha,i',j'|X)}(g(x))$ in the state $(x_{i'},x_{j'})$, so it guarantees $w$ from state $\alpha'_2$.
% Note that if $w$ is the maximum payoff that $\sigma$ guarantees then it is also the maximum payoff that $\sigma'$ guarantees.
Reciprocally given $\sigma'$ a strategy in $\Gamma'$ that guarantees $w'$ from $\alpha'_2$, then $\sigma'^*$ the strategy in $\Gamma'$ such that $\sigma'^*(\alpha'_2)=\sigma'(\alpha'_2)$ and $\sigma'^*(x_{i},.)=i$ also guarantees $w'$ in $\Gamma'(\alpha'_2)$. The strategy $\sigma$ defined by $\sigma(\alpha)=\sigma'(\alpha'_2)$ guarantees the same payoff in the absorbing game. From a strategic point of view the two games are completely equivalent. \hfill $\Box$

%We denote by $\sigma$ the strategy which follows $\sigma'$ as if the state has always been $(\alpha',\alpha')$. For all $h=(i_1,j_1,....,i_n,j_n)\in H_n$, $\sigma(h)=\sigma'( ((\alpha',\alpha'),i_1,j_1,(\alpha',\alpha'),i_2,j_2,...,i_n,j_n,(\alpha',\alpha'))$. Then $\sigma$ guarantees the same payoff as $\sigma'^*$. Therefore it guarantees also $w$.

\subsection{Proof of Theorem \ref{theo2}}\label{preuvecomm}
In this section we prove Theorem \ref{theo2}. Let
$\Gamma=(X,I,J,q,g)$ be a stochastic game where $X$ is a compact
subset of $\R^m$, $I$ and $J$ are finite sets, $q$ is commutative,
deterministic, and $1$-Lipschitz for $\|.\|_1$, and $g$ is
continuous. We will prove that for all $z_1\in \Delta_f(X)$, the
stochastic game $\Gamma(z_1)$ has a uniform value. It is sufficient to
prove that for all $x_1\in X$, $\Gamma(x_1)$ has a uniform value.\\

The outline of the proof is the following. For each $x \in X$ we
separate the action pairs into two different sets. An action pair $(i,j)\in I\times J$ is \emph{cyclic} at $x$ if the play that is obtained by repeating $(i,j)$ starting from $x$, comes back to $x$ after a
finite number of stages. If $(i,j)$ does not satisfy this property,
we say that it is \emph{non-cyclic}.
%For convenience we will speak of cyclic actions and non-cyclic actions in the following even if the property is on the couples.
Denote by $\mathcal{C}(x)$  the set
of cyclic action pairs at $x$ and by $\mathcal{N}\mathcal{C}(x)=(I \times J) \backslash \mathcal{C}(x)$
the set of non-cyclic action pairs at $x$.

We denote by $\Phi_k=\{x ; |\mathcal{C}(x)| \geq k\}$ the set of states with more than $k$ cyclic
action pairs. We will prove by induction on the number of cyclic action pairs that the uniform value exists for all
initial points $x_1\in X$.

We first argue that  $\Phi_{|I \times J|}$ is non-empty and every state $x_1\in \Phi_{|I \times J|}$ has a uniform value. To this end we will note that whatever the players play, only finitely many states can be reached from $x_1$, so that $\Gamma(x_1)$ is in essence a game with a finite number of states. By Mertens and Neyman \cite{Mertens_81} the game has a uniform value.

For the induction step, given a state $x_1$ with $k-1$ cyclic action pairs, we study a family $(\dot{\Gamma}(\varepsilon,x_1))_{\varepsilon>0}$ of games,
which approximate $\Gamma(x_1)$ more and more precisely, and that
have a uniform value. Assume by induction that for all states $x$ in $\Phi_l$, for $l\geq k$, the game $\Gamma(x)$ has a uniform value. For each $\varepsilon>0$, let $\eta$ be defined by uniform continuity of $g$. The game
$\dot{\Gamma}(\varepsilon,x_1)$ is defined as follows: every state $x$ such that there exists $l\geq k$ and $x'\in \Phi_l$ with $\|x-x'\|_1 \leq \eta$ is turned into an absorbing state with payoff the uniform value at $x'$. We will show that $\dot{\Gamma}(\varepsilon,x_1)$ can
be written with a finite number of states, depending on $x_1$.
By Mertens and Neyman \cite{Mertens_81}, it has a uniform value at the initial state $x_1$ denoted $v(\varepsilon)(x_1)$. Finally, we prove that $v(\varepsilon)(x_1)$ converges when $\varepsilon$ goes to $0$ and that the limit is the uniform value of $\Gamma(x_1)$.\newline

We now turn to the formal proof. We first prove an auxiliary Lemma studying the play induced by iterating the same action pair in Section \ref{iterative}. In Section \ref{iniinduction}, we focus on the initial step of the induction. In Section \ref{induction}, we prove the inductive step and conclude the proof.\\

Denote by $q_{i,j}$ the operator from $X$ to $X$ defined by $q_{i,j}(x)=q(x,i,j)$. The map
$q$ is deterministic, so we can define the play along a sequence of
actions. Fix $n \geq 1$ and $h=(i_1,j_1,...,i_n,j_n) \in (I \times
J)^n$. For all integers $l \leq n$ set $x_{l+1}(h)=
q_{i_l,j_l}...q_{i_1,j_1} x_1=\prod_{t=1}^l q_{i_t,j_t} x_1$. We
say that $x$ is \emph{reachable} from $x_1$ if there exists a play from
$x_1$ to $x$.

\subsubsection{Asymptotic behavior of the play induced by one action pair}\label{iterative}

Let $x\in X$. If the action pair $(i,j)$ is cyclic at $x$ then the sequence of states induced by repeating $(i,j)$ from $x$ is periodic. We focus on a non-cyclic action pair at $x$ and we will prove that the set of states along the play induced by repeating $(i,j)$ converges to a periodic orbit of states with strictly more cyclic action pairs than $x$. In order to prove this result, we use the following lemma (Sine \cite{Sine_90}).

\begin{lemma}\label{lemma51}
Let $m\geq 1$, there exists $f(m)\geq 1$ such that for all maps $M$ from $X \subset \R^m$ to $X$, $1$-Lipschtiz for $\|.\|_1$, there exists an integer $L \leq f(m)$ and a family of maps $B_0$,$\cdots$, $B_{L-1}$ such
that
\[
\forall l \in \{0,...,L-1\},\  \lim_{t \rightarrow +\infty} M^{t L + l}=B_l.
\]
\end{lemma}

\noindent A classic example is the case where $M$ is the transition of a Markov chain on a finite set. If $\lambda$ is a complex eigenvalue of $M$ then $|\lambda|\leq 1$ since the map is $1$-Lipschitz. Moreover the theorem of Perron-Frobenius ensures that if $|\lambda|=1$ then there exists $l\leq m$ such that $\lambda^l=1$. The integer $L$ is then the smallest common multiple of all such $l$ and we can take $f(m)=m!$. \newline

%For each $x \in X$ we separate the couples of actions in two different groups.
%On one hand the actions which from $x$ come back to $x$ in less than $\phi(m)$ stages and on the other hand the rest.
%
%\begin{definition}
%Let $x \in X$, the couple $(i,j)\in I \times J$ is cyclic
%in $x$ if there exists $t \leq \phi(m)$ such that $q^t_{i,j} x =x$. We denote by $S(x)$
%the set of cyclic actions in $x$ and by $\overline{S}(x)$ its complementary.
%\end{definition}

%The number of cyclic actions is increasing for the inclusion order along a trajectory and we can prove our result by induction.
%
%\begin{lemma}\label{lemma52}
%If $x \in X$ and $x' \in \Lambda^+_{I \times J}(x,1)$ then $S(x)$, the set of cyclic actions in $x$, is included in $S(x')$, the set of cyclic actions in $x'$.
%\end{lemma}
%
%Proof : Indeed let $x' \in \Lambda^+_{I \times J}(p,1)$, there exists a sequence $(i_1,j_1,...,i_n,j_n) \in (I \times J)^n$ such that $x'=
% \prod_{l=1..n} q_{i_l,j_l} x$. Let $(i^*,j^*) \in S(x)$ and $d \in \N$ such that $q^d_{i^*,j^*}x=x$ then
%
%\[
%q^d_{i^*,j^*} x'  =\prod_{l=1..n} q_{i_l,j_l}  q^d_{i^*,j^*} x
%=x'.
%\]
%%%\end{proofof}

%\begin{example}
%Let $X=\Delta(\X/2\X)$,
%$x_1=(1,0)$, $I=\{1\}$, $J=\{1\}$ and $A=A(1,1)=\begin{pmatrix}
%1/2 & 1/2 \\
%1/2 & 1/2 \end{pmatrix}$. Then $S(x_1)$ is empty and $S(A x_1)=\{(1,1)\}$.
%\end{example}

Applied to our framework, we deduce that, by iterating a non cyclic
action pair $(i,j)$ from $x$, the induced play has a finite number of
limit points with strictly more cyclic action pairs than $x$.

% $(i,j)$ is cyclic at the limit points and, by the
%commutativity assumption, previous cyclic action pairs are still
%cyclic at the limit.
%get that the play, obtained by repeating a non-cyclic couple of actions $(i,j)$ from $x$, converges to a periodic orbit of states with more cyclic actions. Precisely, the ergodic theorem implies that $(i,j)$ becomes cyclic and the commutation assumption ensures that the cyclic actions at $x$ are still cyclic.

\begin{lemma}\label{lemma55}
Let $x\in X$, $(i,j) \in \mathcal{N}\mathcal{C}(x)$ be a non-cyclic action pair at $x$, and $\varepsilon >
0$. There exist an integer $n$ and a finite set $S_x \subset X$ such that
\[
 \forall t\geq n, \exists x' \in S_x, \
\|q_{i,j}^t x -x'\|_1 \leq \varepsilon \text{ and } \sharp \mathcal{C}(x') \geq \sharp \mathcal{C}(x)+1.
\]
\end{lemma}

\noindent{\bf Proof:}
Let $x \in X$, $(i,j)\in
\mathcal{N}\mathcal{C}(x)$ be a non cyclic action pair and $\varepsilon$ be a positive real. We show three properties: first the sequence $(q_{i,j}^t x)_{t\geq 1}$ has a finite number of limit points, then a cyclic action pair at $x$ is still cyclic at the limit points and finally the pair $(i,j)$ becomes cyclic at the limit points. Therefore, the number of cyclic action pairs strictly increases. \\

By Lemma \ref{lemma51} applied to $Q=q_{i,j}$, there exist an integer $L$ and some operators $B_0$,..., $B_{L-1}$ such that
\[
\forall l \in \{0,...,L-1\}\ \lim_{t \rightarrow +\infty} Q^{t L
+ l}=B_l.
\]
In addition, for every $l \in\{0,...,L-1\}$, $Q^l B_0= B_0 Q^l =B_l.$ By compactness of $X$, $B_0 x$ is in $X$ and  there exists an integer $n$ such that
\[
\forall t \geq n, \|Q^{t L} x-B_0 x \|_1 \leq \varepsilon.
\]
Since $Q$ is $1$-Lipschitz for the norm $1$, $\|Q^{t L+l} x- B_l x \|_1 \leq \varepsilon$. Denoting $n'=n(L+1)$ and $S_x=\{B_l x,\ l=0,\ldots,L-1\}$, we have
\[
\forall t\geq n', \ \exists x' \in S_x ,\  \|Q^t x
-x'\|_1 \leq \varepsilon.
\]
The play has a finite number of limit points.\\

Let $(i',j')$ be a cyclic action pair in $x$ and
$d$ an integer such that $q_{i',j'}^d x=x$. We check that $(i',j')$ is still cyclic at the limit points. For all $l\in \{0,...,L-1\}$, we have
\begin{align*}
q_{i',j'}^d B_l x&= \lim_t q_{i',j'}^d Q^{t L+l} x   \\
 &= \lim_t Q^{tL+ l} q_{i',j'}^d x= \lim_t Q^{tL+l} x=B_l x.
\end{align*}
The commutation assumption implies the second equality. Therefore $(i',j')$ is still a cyclic action pair on the set $S_x$.\\
%\begin{align*}
%y A(i',j')^d =x B_0 q_{i',j'}^d &= \lim_n x M^{n L} A(i',j')^d \\
% &= \lim_n x A(i',j')^d M^{n L}= \lim_n x M^{n L}=y
%\end{align*}

Moreover the iterated action pair $(i,j)$, which was non-cyclic at $x$, becomes cyclic at $x'$ for all $x'\in S_x$. For all $l\in \{0,...,L-1\}$, we have
\[
 Q^L B_l x = \lim_t Q^L Q^{t L+l} x = \lim_t Q^{(t+1)L+l} x= B_l x.
\]
All cyclic action pairs at $x$ are still cyclic on $S_x$ and $(i,j)$ becomes cyclic, so the number of cycling action pairs is strictly increasing.
% is in $\mathcal{C}(y)$ but not in $\mathcal{C}(x)$ so $\sharp \mathcal{N}\mathcal{C}(y) < \sharp \mathcal{N}\mathcal{C}(x)$ and since $S$ is increasing along the trajectory we proved that for all $l \in \{0,...,L-1 \}$, $\sharp \mathcal{N}\mathcal{C}(y M^l)  < \sharp \mathcal{N}\mathcal{C}(x)$. There are less non-cyclic actions than in $x$.
\hfill $\Box$

\begin{example}
\rm Consider a stochastic game with state space $X=\Delta(\Z/2\Z)$, initial state $x_1=(1,0)$, trivial sets of actions $I=\{i_1\}$, $J=\{j_1\}$, and transition
\[
Q=q_{i_1,j_1}=\begin{pmatrix}
1/4 & 3/4 \\
3/4 & 1/4 \end{pmatrix}.\]
Then for all $t \in \N$, $Q^t x_1$ has no cyclic action pairs but it converges to $x_{\infty}=(1/2,1/2)$ where the action pair  $(i_1,j_1)$ is cyclic.
\end{example}

% Let $\Gamma=(X,I,J,q,g)$ be a stochastic
%game such that $X$ is a compact subset of $\R^m$, $I$ and $J$ are
%finite sets, $q$ is commutative, deterministic and $1$-Lipschitz for
%$\|.\|_1$, and $g$ is continuous.

%and we say that $x'$ is between $x_1$ and $x$ if there exists a play between $x_1$ and $x$ passing in state $x'$ in between.
%We denote by $M$ an upper bound on the number of stages necessary for cyclic actions to comes back.
\subsubsection{Initialization of the induction}\label{iniinduction}

%In this section, we prove the existence of the uniform value in every state by induction on the number of cyclic action pairs. We first prove the initialization.
\begin{proposition}
The set $\Phi_{|I \times J|}$ is non-empty.
\end{proposition}

The proposition is an immediate corrolary of Lemma \ref{lemma55}. Starting from
any initial state $x^1\in X$, we apply Lemma \ref{lemma55} to one
non-cyclic action pair and we get a state $x^2 \in X$ with
more cyclic action pairs. Then, we can repeat from this new state
and iterate the lemma until all the action pairs are cyclic.

\begin{proposition}\label{prop51}
$\forall x_1\in \Phi_{|I \times J|}$, the game $\Gamma(x_1)$ has a uniform value.
\end{proposition}

%\begin{lemma}\label{lemma53}
%The number of stages needed to reach all points reached by
%iteration of actions in $S(x)$ is finite.
%\end{lemma}

\noindent{\bf Proof:}
Fix $x_1\in \Phi_{|I \times J|}$. Let $M\geq 1$ be such that for all action pairs $(i,j)$, the play that starts at $x_1$ and in which the players repeatedly play $(i,j)$ returns to $x$ after at most $M$ stages. We argue by contradiction that all states reachable from
$x_1$ can be reached in less than $(M-1)\sharp (I \times J)$ stages.
%\[
%\Lambda^+_{S(x)}(x,1)=\Lambda^-_{S(x)}(x, (\phi(m)-1)\sharp
%S(x))
%\]
%By definition the set of points reached in less than $(\phi(m)-1)\sharp S(x)$ is included in the set of points reached without limitation on the number of stages. We how the other inclusion by
%contradiction.
By contradiction let $x^*$ be a state, which is not reached in
$(M-1) \sharp (I \times J)$ stages. We define
\[
t^*=\inf_{t \geq 1}\left \{ t,\ \exists h=(i_l,j_l)_{l=1...t} \in (I \times J)^t,\
x_t(h)=x^* \right \}
\]
the minimum number of stages needed to reach $x^*$. By assumption,
$t^*>(M-1) \sharp (I \times J)$ and
\begin{align*}
 \sum_{(i,j) \in \mathcal{C}(x_1)} \sharp\{l, (i_l,j_l)=(i,j) \} & =t^*  \\
 \Rightarrow \exists (i^*,j^*) \in \mathcal{C}(x_1) \ \sharp\{l,(i_l,j_l) =(i^*,j^*)\} & \geq \frac{t^*}{\sharp
(I \times J)} \\
 \Rightarrow \exists (i^*,j^*) \in \mathcal{C}(x_1) \ \sharp \{l,(i_l,j_l)=(i^*,j^*)\}  & \geq  M.
\end{align*}

So one action pair is repeated more than $M$ times. By
definition, there exists $d^* \leq M$ such that $q^{d^*}_{i,j}
x_1=x_1$. Hence the state at stage $t^*-d^*$ along the sequence of
actions deduced from $h$, by deleting $d^*$ times the action pairs $(i^*,j^*)$, is $x^*$. This contradicts the definition of
$t^*$. Therefore, all states are reached in less than $(M-1) \sharp
(I \times J)$ stages and since $I$ and $J$ are finite, the game
$\Gamma(x_1)$ can be defined only with a finite number of states.

Formally, the game $\Gamma(x_1)$ is a stochastic game with a finite
set of states and finite sets of actions, so it has a uniform value
by the theorem of Mertens and Neyman \cite{Mertens_81}. \hfill
$\Box$
\newline

%\begin{remark}
% Note that even if for every $x\in X$, $\Gamma(x)$ can be defined only with a finite number of state. This set of state depends on $x$. Given $x$ and $x'$ distinct, the minimal set of state necessary to describe the game changes.
%\end{remark}
%

\subsubsection{ Inductive step}\label{induction}

We now prove the inductive step. Fix $0 < k \leq |I \times J|$ and assume that for all $x\in \cup_{l=k}^{|I \times J|} \Phi_l$, the game $\Gamma(x)$ has a uniform value. Fix $x_1 \in \Phi_{k-1}$. \\

%the set of state at more than $\eta$ from the state with more than $k$ cyclic actions. Let $\eta$ be a positive number, for each $x\in X$, we define $\Xi^{\eta}(x)$ the set of states with more than $k$ cyclic actions in the $\eta$-neighbourhood of $x$:
%\[\Xi^{\eta}(x)=\{x' \in X \text{ such that } \sharp \mathcal{C}(x') \geq k \text{ and } \|x -x'\|_1 \leq \eta \}.\]
%By the induction assumption, we know that games starting from these states have a uniform value. Let $\Phi$ be the set of points $x$ reachable from $x_1$ and where $\Xi^{\eta}(x)$ is empty.
First, we check that the $1$-Lipschitz transition and the uniform continuity
of the payoff imply the continuity of the payoff that a
player can guarantee, then we describe the family of auxiliary games
and conclude the proof.

\begin{lemma}\label{lemma54}
Given $\varepsilon>0$, there exists $\eta>0$ such that
if $x\in X$ and player $1$ guarantees $w$ in $\Gamma(x)$ then, for all $x'$, such
that $\|x-x'\|_1\leq \eta$, he guarantees $w-\varepsilon$ in
$\Gamma(x')$.
\end{lemma}

\noindent{\bf Proof:} Given $\varepsilon>0$, for all $(i,j)\in
I\times J$, the map $g(\cdot,i,j)$ is uniformly continuous. Moreover,
the number of maps is finite, so there exists $\eta>0$ such that for
all $x,x'\in X$ with $\|x-x'\|_1 \leq \eta$, we have
\[
\forall (i,j) \in (I \times J), \ |g(x,i,j)-g(x',i,j)| \leq \varepsilon.
\]

We first check the result for pure strategies. Fix $x\in X$. Let $\sigma\in \Sigma$ be a pure strategy, we define $\widetilde{\sigma}(x)$ to be the strategy which plays as if the game were $\Gamma(x)$ no matter what the initial state is. In particular, this strategy does not depend on the state and only on the sequence of actions. % such that for all $\tau \in \T$ the probability on the histories under $(x,\sigma,\tau)$ and $(x,\sigma^*,\tau)$ are the same.
Let $\tau\in J^\N$ be a sequence of actions of player $2$.

We denote by $x_t$ the state at stage $t$ along $(x,\sigma,\tau)$
and $x'_t$ the state at stage $t$ along $(x',\widetilde{\sigma}(x),\tau)$. For
all $(i,j)\in I\times J$, $q$ is a $1$-Lipschtiz function so for all
$t\geq 1$, $\|x_t-x'_t\|_1 \leq \|x-x'\|_1\leq \eta$, and for all $n\geq 1$,
\begin{align*}
|\gamma_n(x,\sigma,\tau)-\gamma_n(x',\widetilde{\sigma}(x),\tau)| & \leq \frac{1}{n}\sum_{t=1}^n |g(x_t,i_t,j_t)-g(x'_t,i_t,j_t)| \\
& \leq \varepsilon.
\end{align*}

Let $\sigma^*$ be a mixed strategy\footnote{Recall that by Kuhn's theorem, a behavioral strategy is equivalent to a mixed strategy.}. We denote by $\PP_{\sigma^*}$ the probability distribution induced by $\sigma^*$ on the set of pure strategies and $\E_{\sigma^*}$ the corresponding expectation. We define the mixed strategy $\widetilde{\sigma^*}$ by associating to each
pure strategy $\sigma$ the strategy $\widetilde{\sigma}(x)$. It is measurable and we have
\begin{align*}
|\gamma_n(x,\sigma^*,\tau)-\gamma_n(x',\widetilde{\sigma^*},\tau)| & \leq \left| \E_{\sigma^*}\left(\gamma_n(x,\sigma,\tau)-\gamma_n(x',\widetilde{\sigma}(x),\tau) \right) \right|\\
& \leq \E_{\sigma^*}\left( \left| \gamma_n(x,\sigma,\tau)-\gamma_n(x',\widetilde{\sigma}(x),\tau) \right| \right) \\
& \leq \varepsilon .
\end{align*}
If player $1$ guarantees $w$ in $\Gamma(x)$ then he guarantees $w-\varepsilon$ in the game $\Gamma(x')$ for every $x'$ such that $\|x-x'\|_1 \leq \eta$. \hfill $\Box$\\

Let $\varepsilon$ be a positive real and $\eta$ be associated to $\varepsilon$ by Lemma \ref{lemma54}.
We denote by $\Phi(\eta)$ the set of states reachable from $x_1$ such that there is no state $x\in \cup_{l=k}^{|I \times J|} \Phi_l$ in the $\eta$-neighbourhood,
\[\Phi(\eta)=\{x \text{ reachable from }x_1,\
\forall x'\in X \ x' \notin \cup_{l=k}^{|I \times J|} \Phi_l \text{ or } \|x -x'\|_1 > \eta \}.\]

%We now prove that $\Phi(\eta)$ is finite.\newline

\begin{proposition}\label{prop52}
The set $\Phi(\eta)$ is finite.
\end{proposition}

\noindent{\bf Proof:} We first prove that there exists $M$ such that any state in $\Phi(\eta)$ can be reached in less than $M$ stages and then we deduce that $\Phi(\eta)$ is finite.\\

%this set of states can be reached in a bounded number of stages. Let
%$H=\{ h \in (I \times J)^{\N} |\ \exists t \geq 1, x_t(h) \in \Phi
%\}$ be the set of possible histories associated to states in $\Phi$.

For each action pair $(i,j)$ in $\mathcal{N}\mathcal{C}(x_1)$,
we denote by $u(i,j)$ the integer given by Lemma \ref{lemma55}. Since there is a finite number of action pairs, there exists $M'$ an integer such that for all $(i,j)\in
\mathcal{N}\mathcal{C}(x_1)$, $u(i,j)\leq M'$ and for all $(i,j)\in
\mathcal{C}(x_1)$, the minimal period of $(i,j)$ is smaller than
$M'$. Set $M=M' \sharp (I \times J)$.

%For all $x\in \Phi(\eta)$, we denote by $t^*(x)=\inf \{t |\ \exists h \in (I
%\times J)^t \ x_t(h)=x \}$, the least number of stages necessary to
%reach $x$.

We prove that for all $x\in \Phi(\eta)$, $t^*(x)=\inf \{t |\ \exists h \in (I
\times J)^t \ x_t(h)=x \}$, the least number of stages necessary to
reach $x$, is smaller than $M$.\\

By contradiction, let $x\in \Phi(\eta)$ such that $t^*=t^*(x) \geq M$ and $h$ be an
history associated to $x$ and $t^*$, then one action pair $(i^*,j^*)$ is
repeated more than $M'$ times. This action pair is either cyclic or non-cyclic at $x_1$. If this action pair is cyclic, the history can be shortened, as
in the proof of Proposition \ref{prop51}, which is absurd with
respect to the definition of $t^*$. If this action pair is non-cyclic at $x_1$,
there exists $\overline{x} \in X$ such that
\begin{align*}
             & \|q_{i^*,j^*}^{M'} x_1 -\overline{x} \|_1 \leq \varepsilon,      \\
\text{ and } & \sharp \mathcal{C}(\overline{x}) > k-1.
\end{align*}
Denote by $h'$ the sequence of action pairs
where $(i^*,j^*)$ has been deleted $M'$ times from $h$ and $x'$
the state obtained from $\overline{x}$ by playing $h'$. The transition is $1$-Lipschitz and $\mathcal{C}$ is non-decreasing, therefore we have
\begin{align*}
             & \|x-x'\|_1 \leq \varepsilon,      \\
\text{ and } & \sharp \mathcal{C}(x') > k-1,
\end{align*}
which contradicts the definition of $x.$\\

To conclude notice that there exists a finite number of actions, therefore the set $\Phi(\eta)$ is finite. \hfill $\Box$ \\

%
%
%We now define the auxiliary game by choosing for each
%$\varepsilon>0$, an $\eta$-neighbourhood small enough.
%maand then we show that a player can guarantee a uniformly continuous payoff and we use this continuity to forget what is happening in the neighbourhood of $\phi_k$ and et a game with a finite number of states.

%, so for all $x'\in X$ and any $\varepsilon$-optimal strategy of $\Gamma(x')$, this strategy in the games beginning in the neighbourhood of $x'$.

By Proposition \ref{prop52}, the set of states, reachable from $x_1$,
 and at a distance at least $\eta$ from any state with more than $k$ cyclic
action pairs, i.e. $\Phi(\eta)$, is finite. We denote by
$q(\Phi(\eta))$ the set of all states obtained by one transition
from one of these states and, which are not already in $\Phi(\eta)$.
The set $q(\Phi(\eta))$ is finite and for each $x\in q(\Phi(\eta))$,
there exists $\xi(x) \in \cup_{l=k}^{|I \times J|} \Phi_l$ such that $d(x,\xi(x))\leq \eta$. The induction assumption implies therefore that the game
$\Gamma(\xi(x))$ has a uniform value denoted by $v^*(\xi(x))$. We define the auxiliary game
$\dot{\Gamma}(\varepsilon,x_1)$ as follows: the initial state is
$x_1$, the sets of actions are $I$ and $J$, the transition function and
reward functions are given by:
\begin{small}
\begin{center}
\begin{tabular}{rl}
 $\dot{q}(x,i,j)=$ & $\begin{cases}
                       q_{i,j} x & \text{ if }x\in \Phi(\eta)\\
                       x & \text{ if }x \in q(\Phi(\eta)) \\
                       x & \text{ otherwise,}\\
                     \end{cases}$\\
and $\dot{g}(x,i,j)=$ & $\begin{cases}
                       g(x,i,j) & \text{if }x\in \Phi(\eta)\\
                       v^*(\xi(x)) & \text{ if }x \in q(\Phi(\eta)) \\
                       0 & \text{ otherwise}.
                       \end{cases}$
\end{tabular}
\end{center}
\end{small}
The sets of strategies for players $1$ and $2$ are the same as in the game $\Gamma$. In the game starting at $x_1$, all the states are in $\Phi(\eta)$ or $q(\Phi(\eta))$. Since both sets are finite, this game is formally a stochastic game with a finite set of states and finite sets of actions. Therefore $\dot{\Gamma}(\varepsilon,x_1)$ has a uniform value by the theorem of Mertens and Neyman \cite{Mertens_81}.

\begin{proposition}\label{prop53}
 $\dot{\Gamma}(\varepsilon,x_1)$ has a uniform value in $x_1$ denoted by $v^*(\varepsilon)(x_1)$.
\end{proposition}

We now prove that when $\varepsilon$ goes to $0$, the value $v^*(\varepsilon)(x_1)$ has to converge and the limit is the uniform value of the game $\Gamma(x_1)$. We first prove that the value of the auxiliary game is a good approximation to what the players can guarantee in $\Gamma(x_1)$.
\begin{proposition}\label{prop54}
  If player $1$ can guarantee $w$ in
  $\dot{\Gamma}(\varepsilon,x_1)$ then he can guarantee $w- 3
  \varepsilon$ in $\Gamma(x_1)$.
\end{proposition}

\noindent{\bf Proof:}
By assumption, there exists $\dot{\sigma}$  a strategy of player $1$
in $\dot{\Gamma}(\varepsilon,x_1)$ and a stage $\dot{N}\geq 1$ such that
\[
\forall n \geq \dot{N}, \ \forall \dot{\tau}, \
\dot{\gamma}_n(x_1,\dot{\sigma},\dot{\tau}) \geq w- \varepsilon.
\]
For each state $x\in q(\Phi(\eta))$, we denote by $\sigma^{\xi,x}$
the strategy given by Lemma \ref{lemma54} with respect to the
point $\xi(x)$ and to an $\varepsilon$-optimal strategy in
$\Gamma(\xi(x))$ such that
\[
\exists N(x)\geq 1, \ \forall n \geq N(x), \  \forall
\tau, \  \gamma_n(x,\sigma^{\xi,x},\tau) \geq v^*(\xi(x))-2 \varepsilon.
\]
Let $\overline{N}=\max(N(x),x\in \Phi(\eta))$ be an upper bound.\\

Given an infinite play $h \in (X\times I \times J)^{\N}$, we denote by $\theta(h)$ the first stage where the state is at a distance less than $\eta$ from a state in $\cup_{l=k}^{|I \times J|} \Phi_l:$
\[\theta(h)=\inf_{t \geq 1}\{t | x_t(h) \in q(\Phi(\eta)) \}.\]
We define the strategy $\sigma$ which plays optimally in $\dot{\Gamma}$ until a state $x'\in q(\Phi(\eta))$ is reached, and then optimally as if the remaining game was $\Gamma(\xi(x'))$. Formally, we have
\begin{align*}
\forall n\geq 1, \  \sigma_n(h) = \begin{cases}
                  \dot{\sigma}_n(h) & \text{ if } n \leq \theta(h)-1 \\
                  \sigma_{n-\theta(h)+1}^{\xi,x_{\theta(h)}(h)} & \text{ if } n \geq \theta(h).
               \end{cases}
\end{align*}

\noindent We prove that $\sigma$ guarantees $w-3\varepsilon$. Let $\tau$ be a strategy of player $2$, we denote by $x_t$ the state at stage $t$. Let $N^*\in \N$ such that $N^* \geq \dot{N}$ and $\frac{\overline{N}}{N^*} \leq \varepsilon$. Fix $n \geq N^*$, we separate the histories in two groups depending on whether $n-\theta(h)+1>\overline{N}$ or $n-\theta(h)+1 \leq \overline{N}$.\\

We first focus on the set of histories  $\{ h \in H_\infty,\ n-\theta(h)+1>\overline{N}\}$ and notice that on these histories the expected average payoff between $\theta(h)$ and $n$ is close to the uniform value at $\xi(x_\theta(h))$.

We denote by $\sigma^{h_n}$ and $\tau^{h_n}$ the strategies induced by $\sigma$ and $\tau$ given the finite history $h_n$. Since $\|x_{\theta(h)}-\xi(x_{\theta(h)})\| \leq \eta $, we have
\begin{align*}
&\E_{x_1,\sigma,\tau} \left( \sum_{t=\theta(h)}^n g(x_t,i_t,j_t) \1_{n-\theta(h)+1>\overline{N}} \right) \\
=&  \E_{x_1,\sigma,\tau} \left( \gamma_{n-\theta(h)+1}(x_{\theta(h)},\sigma^{h_{\theta(h)}},\tau^{h_{\theta(h)}})(n-\theta(h)+1) \1_{n-\theta(h)+1>\overline{N}} \right) \\
\geq &  \E_{x_1,\sigma,\tau} \left( \left( v^*(\xi(x_{\theta(h)}))-2 \varepsilon \right)(n-\theta(h)+1) \1_{n-\theta(h)+1>\overline{N}} \right).
\end{align*}
Therefore
\begin{align*}
& \frac{1}{n}  \E_{x_1,\sigma,\tau} \left( \sum_{t=1}^n g(x_t,i_t,j_t) \1_{n-\theta(h)+1>\overline{N}} \right)\\
& = \frac{1}{n}  \E_{x_1,\sigma,\tau} \left(\left( \sum_{t=1}^{\theta(h)-1} g(x_t,i_t,j_t)+ \sum_{t=\theta(h)}^n g(x_t,i_t,j_t) \right) \1_{n-\theta(h)+1>\overline{N}} \right) \\
& \geq  \E_{x_1,\sigma,\tau} \left( \frac{1}{n} \left(
  \sum_{t=1}^{\theta(h)-1}
  g(x_t,i_t,j_t)+ v^*( \xi(x_{\theta(h)}))(n- \theta(h)+1) \right) \1_{n-\theta(h)+1 \geq \overline{N}}- 2 \varepsilon \1_{n-\theta(h)+1 \geq \overline{N}} \right).
\end{align*}
%
%Let $h$ be an infinite history such that $n-\theta(h)>\overline{N}$ then for each history there are two cases. On one hand if $n-\theta(h)>\overline{N}$, $\sigma$ has played optimally in the game from the state $\xi(x_\theta)$ for long enough in order for the payoff to be above $v^*(\xi(x_\theta))-\varepsilon$. On the other hand if $n-\theta(h)<\overline{N}$, then the part of the play after $\theta(h)$ has a weight  less than $\varepsilon$ of the total. We split the payoff depending on this criteria,
%%
%\begin{align*}
%\gamma_N(x,\sigma,\tau) & =\frac{1}{N}  E_{x,\sigma,\tau} \left( \psi(x_n,i, j) \1_{N-\tht \geq \overline{N}} + \psi(x_n, i, j) \1_{N-\tht <\overline{N}} \right), %\\
%\end{align*}
%
%
%
%& =\frac{1}{N}  E_{x,\sigma,\tau} \left(
%\sum_{n=1}^{\tht} r(x_n,i_n,j_n)+ \sum_{n=\tht+1}^N r(x_n,i_n,j_n) \right) \\
%%& =\frac{1}{N}  E_{x,\sigma,\tau} \left( \psi(x_n,
%%i, j) \right) \\
%\begin{align*}
%\gamma_N(x,\sigma,\tau)
%& =\frac{1}{N}  E_{x,\sigma,\tau} \left(
%\sum_{n=1}^{\tht} r(x_n,i_n,j_n)+ \sum_{n=\tht+1}^N r(x_n,i_n,j_n) \right) \\
%& =\frac{1}{N}  E_{x,\sigma,\tau} \left( \psi(x_n,
%i, j) \right) \\
%\end{align*}
%We study both parts separately and focus first
%

\noindent We now consider the set of histories $\{ h \in H_\infty,\ n-\theta(h)+1 \leq \overline{N}\}$ and notice that the payoff between $\theta(h)$ and $n$ has a small weight.  By definition on this set of histories
\[
\frac{n-\theta(h)+1}{n}\leq \frac{\overline{N}}{N^*} \leq \varepsilon.
\]
Moreover we have
\[
\forall x\in X,\ \forall x'\in q(\Phi(\eta)),\ \forall i\in I,\ \forall j\in J,\ g(x,i,j)\geq -1 \geq v(\xi(x'))-2.
\]
It follows that
\begin{align*}
& \frac{1}{n} \E_{x_1,\sigma,\tau} \left( \sum_{t=1}^n g(x_t,i_t,j_t) \1_{n-\theta(h) \leq \overline{N}} \right) \\
& = \E_{x_1,\sigma,\tau} \left( \frac{1}{n} \left( \sum_{n=1}^{\theta(h)-1} g(x_t,i_t,j_t)+\sum_{t=\theta(h)}^n
g(x_t,i_t,j_t) \right) \1_{n-\theta(h)+1 \leq \overline{N}} \right) \\% ligne 2
% & = E_{x,\sigma,\tau}  \left( \frac{1}{N} \left( \sum_{n=1}^{\theta(h)}
% r(x_n,i_n,j_n)+ v(\xi(\xt_{\theta(h)}))(N- \theta(h))+\sum_{n=\theta(h)+1}^N
%  \left( r(x_n,i_n,j_n)- v( \xi(x_{\theta(h)})) \right) \right) \1_{N-\theta(h) < %\overline{N}} \right) \\ %ligne 3
 & \geq  \E_{x_1,\sigma,\tau} \left(\frac{1}{n} \left(
  \sum_{t=1}^{\theta(h)-1} g(x_t,i_t,j_t) + v^*( \xi(x_{\theta(h)}))(n- \theta(h)+1)
  - 2 (n-\theta(h)+1) \right) \1_{n-\theta(h)+1 \leq \overline{N}} \right) \\ % ligne 4
% & \geq  E_{x,\sigma,\tau} \left(\frac{1}{N} \left(
%  \sum_{n=1}^{\theta(h)} r(x_n,i_n,j_n) + v( \xi(x_{\theta(h)}))(N- \theta(h))
%  - 2 (N-\theta(h)) \right) \1_{N-\theta(h) < \overline{N}} \right) \\ % ligne 4
 & \geq  \E_{x_1,\sigma,\tau} \left( \frac{1}{n} \left(
  \sum_{t=1}^{\theta(h)-1}
  g(x_t,i_t,j_t)+ v^*( \xi(x_{\theta(h)}))(n- \theta(h)+1) \right) \1_{n-\theta(h)+1 \leq \overline{N}}- 2 \varepsilon \1_{n-\theta(h)+1 \leq \overline{N}} \right).
\end{align*}
Therefore by summing the two inequalities we get the result
\[
\gamma_n(x_1,\sigma,\tau) \geq \dot{\gamma}_n(x_1,\dot{\sigma},\tau)
- 2\varepsilon \geq w- 3\varepsilon.
\]
\hfill $\Box$

%To conclude the proof, we show that $v(\varepsilon)(x_1)$, the value
%of the auxiliary game $\dot{\Gamma}(x_1,\varepsilon)$, converges
%when $\varepsilon$ converges to $0$, and the limit is the value of
%the game $\Gamma(x_1)$.

It follows from Proposition \ref{prop54} that for all
$\varepsilon>0$, player $1$ can guarantee
$v(\varepsilon)(x_1)-3\varepsilon$ in the game $\Gamma(x_1)$. So
player $1$ can guarantee the superior limit when $\varepsilon$
converges to $0$: for all $\delta>0$, there exists $n_1$ and a
strategy $\sigma^* \in \Sigma$ such that for all $\tau \in \T$, for
all $n'\geq n_1,$
\[
\gamma_{n'}(x_1,\sigma^*,\tau) \geq \limsup_{\varepsilon \rightarrow 0} v(\varepsilon)(x_1)-\delta.
\]
The same argument shows that player $2$ can guarantee the inferior
limit. Therefore, for all $\delta>0$, there exist $n_2$ and a
strategy $\tau^*\in \T$ such that for all $\sigma \in \Sigma$, for
all $n'\geq n_2,$
\[
\gamma_{n'}(x_1,\sigma,\tau^*) \leq \liminf_{\varepsilon \rightarrow 0} v(\varepsilon)(x_1)+\delta.
\]
Given $\delta>0$ and $n'\geq \max(n_1,n_2)$, we have
\[
\limsup_{\varepsilon \rightarrow 0} v(\varepsilon)(x_1)-\delta \leq \gamma_{n'}(x_1,\sigma^*,\tau^*) \leq \liminf_{\varepsilon \rightarrow 0} v(\varepsilon)(x_1)+\delta.
\]
Therefore $v(\varepsilon)(x_1)$ converges when $\varepsilon$ goes to $0$ and the limit is the uniform value of the game $\Gamma(x_1)$. This proves the induction hypothesis at the next step and concludes the proof. For all $x_1\in X$, the game $\Gamma(x_1)$ has a uniform value.
%and for all initial probability with finite supports, $z_1\in \Delta_f(X)$, $\Gamma(z_1)$ has a uniform value.

\subsection{Proof of Corollary \ref{mdp_noir}}\label{preuveblindstate}

In this section, we provide a short proof of Corollary \ref{mdp_noir}. Recall that given a state-blind repeated game $\Gamma^{sb}=(K,I,J,q,g)$ with a commutative transition $q$, we define the auxiliary stochastic game
$\Psi=(X,I,J,\widetilde{q},\widetilde{g})$ where $X=\Delta(K)$, $\widetilde{q}$ is the linear extension of $q$, and $\widetilde{g}$ is the linear extension of $g$.

%A state, in this new game, is the common
%belief of the players over the state in $\Gamma^{sb}(K,I,J,q,g)$.
%Since $K$ is finite,  $X$ can be embedded in $\R^{K}$ and the transition $\widetilde{q}$ is deterministic, $1$-Lipschitz for $\|.\|_1$, and commutative. Furthermore, $\widetilde{g}$ is continuous and therefore we can apply Theorem \ref{theo2} to $\Psi$. It follows that for each initial state $p_1 \in X$, $\Psi(p_1)$ has a uniform value. We will check that it is  the uniform value of the state-blind repeated game $\Gamma^{sb}(p_1)$ and deduce the following corollary. 

In this framework deducing the existence of the uniform value in the original repeated game  from the existence of the uniform value in the auxiliary game is easy since the sets of strategies are almost the same in the two games. A player can use a strategy of the repeated game $\Gamma$ in $\Psi$ by looking only at the actions played and reciprocally a player can use a strategy of the stochastic game $\Psi$ in the repeated game $\Gamma$ by completing the sequence of actions with the unique sequence of compatible beliefs.\

% with full monitoring which fulfills the assumption of the theorem. In $\widetilde{\Gamma}$, given a couple of actions, the transition $\widetilde{q}$ is deterministic so we can focus on strategies which depend only on the actions. This set of strategies is exactly the set of strategies of $\Gamma^{sb}$. Moreover the transition $\widetilde{q}$ is defined such that the state in $\widetilde{\Gamma}$ is the law of the state in $\Gamma^{sb}$. We deduce that the two games have the same finite values and that the existence of a uniform value in $\widetilde{\Gamma}$ implies the existence of the uniform value in the original game.

%Given a state-blind stochastic game $\Gamma^{sb}(K,I,J,q,r)$ such that $q$ is commutative, we define the auxiliary stochastic game $\widetilde{\Gamma}(p_1)=(X,I,J,\widetilde{q},\widetilde{r},p_1)$ by: $X=\Delta(K)\subset \R^{K}$, $I$ and $J$ are the same as in $\Gamma^{sb}$, $p_1$ is the Dirac mass in $x_1$, $\widetilde{q}$ and $\widetilde{r}$ are the linear extensions of $q$ and $r$ to $X$. Explicitly, let $p\in X$, $i\in I $ and $j\in J$, we have
% \[\widetilde{r}(p,i,j)=\sum_{k\in K} r(k,i,j)\]
%and
%\[
%\widetilde{q}(p,i,j)=\sum_{k\in K} q(k,i,j).\]

\noindent{\bf Proof:} The set of strategies in the game $\Gamma^{sb}$ are respectively
denoted by $\Sigma^{sb}$ and $\T^{sb}$. We will denote in this proof
the payoff in the $n$ stage game by $\gamma_n^{sb}$ and the value of
the $n$-stage game by $v_n^{sb}(p_1)$ for all $n\geq 1$.

We denote by $\widetilde{H_t}$ the set of histories in $\Psi$ of length $t$, by $\widetilde{\Sigma}$ the set of strategies of player $1$, and by $\widetilde{\T}$ the set of strategies of player $2$. Let $p_1 \in \Delta(K)$, $\widetilde{\sigma} \in \widetilde{\Sigma}$ and $\widetilde{\tau} \in \widetilde{\T}$. The payoff in the $n$-stage game, starting from $p_1$ and given that the players follow $\widetilde{\sigma}$ and $\widetilde{\tau}$, is denoted by $\widetilde{\gamma_n}(\delta_{p_1},\widetilde{\sigma},\widetilde{\tau})$ and the value by $w_n(p_1)$. The set $X$ is compact, $\widetilde{g}$ is continuous and the transition $\widetilde{q}$ is commutative and deterministic, so we can apply Theorem \ref{theo2} to $\Psi$. We denote by $w^*(p_1)$ the uniform value. The values of both games coincide since the payoff and strategy sets coincide up to the following identification.\\

%are equal by constructing two bijections between the sets of strategies of player $1$ (resp. player $2$) in both games wich preserve thethey are equal to their equivalent in $\Gamma^{sb}$ by proving there exists functions in-between the two sets of strategies of player $1$ in the two games and in-between the two sets of strategies of player $2$.

We focus on the case of player $1$ since the situation is symmetric for player $2$. Let $\sigma^{sb}$ be a strategy in $\Sigma^{sb}$, then it defines naturally a strategy $\widetilde{\sigma}$ in $\widetilde{\Sigma}$ by forgetting the states. If we denote by $\Pi^t$ the projection from $\widetilde{H_t}$ on $H_t^{sb}$ that keeps only the actions: for all $t\geq 1$, we define
\[\widetilde{\sigma}(\widetilde{h}_t)=\sigma^{sb}(\Pi^t(\widetilde{h}_t)).\]
Reciprocally for all $t\geq 1$, given a sequence of actions $h_t^{sb}=(i_1,j_1,...i_t,j_t)$, the completion $\Xi^t(h^{sb})$ in $\widetilde{H}_t$ is the unique sequence such that $p_1$ is fixed and for all $t\geq 1$, $q(p_t,i_t,j_t)=p_{t+1}$. Let $\widetilde{\sigma}$ be a strategy in $\widetilde{\Sigma}$, then we define the strategy $\sigma^{sb}$ by completing the history: for all $t\geq 1$
\[\sigma^{sb}(h^{sb}_t)=\widetilde{\sigma}(\Xi^t(h^{sb}_t)).\]
\noindent A  similar procedure gives two functions between the sets of strategies of player $2$.\\

Given $\widetilde{\sigma}\in\widetilde{\Sigma}$ and $\tau^{sb}\in \T^{sb}$, set $\sigma^{sb}\in \Sigma^{sb}$ and
$\widetilde{\tau}\in\widetilde{\T}$ as in the previous paragraph. By definition of $\widetilde{q}$, the state at stage $t$ in $\Psi$ under $\P_{\delta_{p},\widetilde{\sigma},\widetilde{\tau}}$ is equal to the law of the state in $\Gamma^{sb}$ under $\P_{p,\sigma^{sb},\tau^{sb}}$. Therefore for all $n\geq 1$, we have
\[
\gamma_n^{sb}(p_1,\sigma^{sb},\tau^{sb})=\widetilde{\gamma_n}(\delta_{p_1},\widetilde{\sigma},\widetilde{\tau}).
\]

Finally, let $\varepsilon>0$, $\widetilde{\sigma}$ be an
$\varepsilon$-optimal strategy in $\Psi$ and $N\geq 1$ an integer
such that for all $\widetilde{\tau} \in \widetilde{\T}$,
\begin{align*}
\widetilde{\gamma_n}(\delta_{p_1},\widetilde{\sigma},\widetilde{\tau}) \geq w^*(p_1)-\varepsilon,
\end{align*}
then for all $\tau^{sb}\in \T^{sb}$, we have
\begin{align*}
\gamma_n^{sb}(p_1,\sigma^{sb},\tau^{sb})& =\widetilde{\gamma_n}(\delta_{p_1},\widetilde{\sigma},\widetilde{\tau})\\
                                                 & \geq w^*(p_1)-\varepsilon.
\end{align*}
The strategy $\sigma^{sb}$ guarantees $w^*(p_1)-\varepsilon$ and therefore player $1$ guarantees $w^*(p_1)$. By symmetry, player $2$ guarantees $w^*(p_1)$ and the game $\Gamma^{sb}(p_1)$ has a uniform value equal to $w^*(p_1)$. \hfill $\Box$

\subsection{Extensions.}\label{extension}

%Finally the theorem is still true if we use a weaker assumption on the transition. We say that  a transition function $q$ weakly commutes on $X$ if for all $x \in X$, for all $i,i' \in I$, there exists $i''\in I$,
%\[
%q(q(x,i'),i)=q(q(x,i),i'').
%\]
%The key point is that given an order on the couple of actions and an initial state, every finite play has an equivalent play which leads to the same state where the couple of actions played are increasing. Especially the occurences of each couple of actions are played consecutively. Thus it is enough to remember how many each actions is played in this equivalent play to recall the transition from the initial state.

The proof of Theorem \ref{theo2} can be extended by replacing some of the lemmas with more general results. The result of Sine \cite{Sine_90}, for example, applies to more general norms than the norm $\|.\|_1$.
\begin{definition}
 A norm on $\R^n$ is polyhedral if the unit ball has a finite number of extreme points.
\end{definition}
For example the norm $\|.\|_1$ and the sup norm are polyhedral norms
but not the Euclidean norm. For polyhedral norm, the application of
the theorem of Sine \cite {Sine_90} to compact sets gives the
following results,
%which is a generalixation of the lemma $1$ and under which our proof still holds.
\begin{lemma}
Let $N(.)$ be a polyhedral norm and $K\subset \R^m$ be a compact
set. There exists $\phi(N,m) \in \N$ such that for all functions $T$,
$1$-Lipschtiz for $N$, there exists $t\leq \phi(N,m)$ such that
$(T^{tn})_{n\in \N}$ converges.
\end{lemma}

We deduce the following theorem.

\begin{theorem}\label{theo3}
Let $\Gamma=(X,I,J,q,g)$ be a stochastic game, such that $X$ is a compact set of $\R^m$, $I$ and $J$ are finite sets, $q$ is commutative deterministic $1$-Lipschitz for a polyhedral norm, and $g$ is continuous. For all $z_1\in \Delta_f(X)$, the stochastic game $\Gamma(z_1)$ has a uniform value.
\end{theorem}

\noindent This theorem does not apply to Example \ref{exemple1} on the circle and the existence of a uniform value in this model is still an open question.

We can obtain new results on non zero-sum stochastic games by replacing the theorem from Mertens and Neyman \cite{Mertens_81} with other existence results. First, Vieille \cite{Vieille_2000a}\cite{Vieille_2000b} proves the existence of an equilibrium payoff in every two-player stochastic games. So our proof, adapted to the non zero-sum case leads to the following result.
\begin{theorem}\label{theo4}
Let $\Gamma=(X,I,J,q,g_1,g_2)$ be a two-player non zero-sum stochastic game such that $X$ is a compact subset of $\R^m$, $I$ and $J$ are finite sets of actions, $q$ is commutative deterministic $1$-Lipschitz for $\|.\|_1$ and $g_1$ and $g_2$ are continuous. Then, for all $z_1\in \Delta_f(X)$, the stochastic game $\Gamma(z_1)$ has an equilibrium payoff.
\end{theorem}

Secondly, there exist some specific classes of $m$-player stochastic games where the existence of an equilibrium has been proven. For example, Flesch, Schoenmakers and Vrieze \cite{Flesch_2008}\cite{Flesch_2009} prove the existence of an equilibrium for $m$-player stochastic games where each player controls a finite Markov chain and the payoffs depend on the $m$ states and the $m$ actions at stage $n$. Note that the commutativity assumption here is reduced to a condition player by player. As in our proof, the commutativity assumption implies that we can study deterministic transitions $1$-Lipschitz for the norm $\|.\|_1$.

\begin{theorem}\label{theo5}
Let $\Gamma=\left((X_j,I_j,q_j)_{j\in \{1,...,m\}},g \right)$ be a $m$-player product-state space stochastic game such that for all $j\in \{1,...,m\}$, $X_j$ is a compact subset of $\R^{m_j}$, $I_j$ is a finite set of actions, $q_j$ is commutative deterministic $1$-Lipschitz for $\|.\|_1$ and $g:\prod(X_j \times I_j)\rightarrow [0,1]^m$ is continuous. For all $z_1\in \Delta_f(\prod_{j} X_j)$, the stochastic game $\Gamma(z_1)$ has an equilibrium payoff.
\end{theorem}

%Finally let us stress out that our proof does not apply to stochastic game with weakly commutative transitions. The number of cyclic actions is not increasing along the trajectories if there are some reset actions as the action $\alpha$ in example \ref{exemple2}. Therefore the induction is not true. It is open if there exists a uniform value in weakly commutative stochastic games where the players do not observe the state.

\section*{Acknowledgments}

I thank J.Flesch, S.Gaubert, J.Renault, E. Solan and S.Sorin as well as three referees. The suggestions
they provided were extremely helpful.

Part of this work has being done in $C \& O$ (University Paris 6) and GREMAQ
(University Toulouse 1 Capitole).

The author gratefully acknowledges the support of the Agence
Nationale de la Recherche, under grant ANR JEUDY, ANR-10-BLAN 0112,
41 as well as the PEPS project Interactions INS2I �Propri\'et\'es des
Jeux Stochastiques de Parit\'e \`a Somme Nulle avec Signaux� and the
Israel Science Foundation, under grant ISF Grant $\# 1517/11$.

\bibliography{biblio_cyclique}

\end{document}